\documentclass[a4paper,12pt]{amsart}

\usepackage{amssymb}
\usepackage{amsmath,amsthm}
\usepackage{enumerate,paralist}
\usepackage{amstext}
\usepackage{dsfont}
\usepackage[english]{babel}
\usepackage{color}
\usepackage{textgreek}

\usepackage[a4paper, left=2.5cm, right=2.5cm, top=3cm, bottom=3cm]{geometry}
\usepackage{mathtools}
\usepackage{latexsym}
\usepackage{mathrsfs}
\usepackage{MnSymbol}
\usepackage{bbm}

\theoremstyle{plain}

\newtheorem{corollary}{Corollary}
\newtheorem{lemma}{Lemma}
\newtheorem{proposition}{Proposition}
\theoremstyle{definition}

\newtheorem{example}{Example}
\newtheorem{remark}{Remark}

\numberwithin{theorem}{section}
\numberwithin{corollary}{section}
\numberwithin{lemma}{section}
\numberwithin{definition}{section}
\numberwithin{example}{section}
\numberwithin{remark}{section}
\numberwithin{proposition}{section}
\numberwithin{assumption}{section}

\usepackage{color}
\definecolor{MY}{rgb}{0.5,0,0.45}

\newcommand{\scr}[1]{\mathscr #1}

\newcommand{\E}{\mathbb{E}}

\newcommand{\1}{\mathds 1}

\newcommand{\eps}{\varepsilon}

\newcommand{\dD}{\mathcal{D}}
\newcommand{\eE}{\mathcal{E}}

\def\al{\alpha}
\def\la{\lambda}
\def\ga{\gamma}

\def\Del{\Delta}

\def\na{\nabla}

\def\eps{\varepsilon}

\def\Om{\Omega}
\def\dOm{\partial \Omega}

\def\pd{\partial}
\def\lo{\lambda_\Omega}

\def\ld{\lambda_{\partial \Omega}}

\def\laa{\lambda_\alpha}

\def\lp{\left(}
\def\rp{\right)}

\def\ff{\frac}
\def\Ric{\mathcal Ric}
\def\rr{\rho}

\begin{document}

\title[]{Functional Inequalities for Brownian motion on manifolds with sticky-reflecting boundary diffusion\footnote{Supported in
 part by  the National Key R\&D Program of China (No. 2022YFA1006000, 2020YFA0712900) and NNSFC (11921001).}}
\author{Marie Bormann}
\address{Universit\"at Leipzig, Fakult\"at f\"ur Mathematik und Informatik, Augustusplatz 10, 04109 Leipzig, Germany and Max Planck Institute for Mathematics in the Sciences, 04103 Leipzig, Germany}
\email{bormann@math.uni-leipzig.de}

\author{Max von Renesse}
\address{Universit\"at Leipzig, Fakult\"at f\"ur Mathematik und Informatik, Augustusplatz 10, 04109 Leipzig, Germany}
\email{renesse@math.uni-leipzig.de}

\author{Feng-Yu Wang}
\address{Center for Applied Mathematics, Tianjin University, Tianjin 300072, China}
\email{wangfy@tju.edu.cn}

\begin{abstract}
We prove geometric upper bounds for the Poincar\'{e} and Logarithmic Sobolev constants for Brownian motion on manifolds with sticky reflecting boundary diffusion i.e.\ extended Wentzell-type boundary condition under general curvature assumptions on the manifold and its boundary. The method is based on an interpolation involving energy interactions between the boundary and the interior of the manifold. As side results we  obtain explicit geometric bounds on the first nontrivial Steklov eigenvalue,  for the norm of the boundary trace operator on Sobolev functions,  and on the boundary trace logarithmic Sobolev constant. The  case of Brownian motion with pure sticky reflection is also treated.
\end{abstract}

\date{\today}

\maketitle

\section{Introduction and description of main results}

Let $\Om$ be a smooth compact connected Riemannian manifold of dimension $d\ge2$ with  smooth connected
boundary $\dOm$. We consider the semigroup  on $C(\Om)$ induced from the  Feller generator $(\dD(A),A)$ given by
\begin{align*}
	\dD(A) &= \{f\in C(\Om)\ |\ Af\in C(\Om)\} \\
	Af &= \Del f\1_\Om + \lp \beta\Del^\tau f - \ga \frac{\pd f}{\pd N}\rp\1_{\dOm},
\end{align*}
where $\frac{\pd f}{\pd N}$ is the outer normal derivative, $\Del^\tau$ is the Laplace-Beltrami operator on $\dOm$, $\beta \geq 0$ and $\ga>0$. The induced Markov process is a diffusion on $\Om$ which performs Brownian motion in the interior and sticky reflected Brownian motion along the boundary. The case of pure sticky reflection but no diffusion along the boundary corresponds to the regime $\beta =0$.  First rigorous constructions of such processes on special domains  $\Omega$ date back to  \cite{MR126883,MR929208,MR287612} and were later extended to jump-diffusion processes  on general domains cf. \cite{MR245085}, \cite{MR4176673}.
 Renewed interest for models with $\beta >0$ emerges from applications in interacting particle systems with singular boundary or zero-range pair interaction \cite{MR4096131,MR2198199,MR2215623,konarovskyi2017reversible,nonnenmacher2018overdamped}.  An efficient construction in the diffusion case via Dirichlet forms was given recently in \cite{MR3498008,gv}, the associated JKO-Wasserstein gradient flow structure is investigated in \cite{leomon}.

\smallskip
In this work we come back to the problem raised in \cite{mvr} of estimating the speed of convergence of this family of processes to their equilibrium state. Any such equilibrium  is  a composition  of the two mutually singular uniform measures on the interior and the boundary and thus has infinite negative curvature even  if $\Om$ and $\dOm$ are positively curved. As a consequence, standard arguments for e.g.\ log-concave measures are not applicable. However,  as shown in \cite{mvr} one can use a simple interpolation together with the Reilly formula to obtain explicit estimates for the $L^2$-spectral gap in the positive curvature case. We aim to extend the interpolation  approach in two ways.
First we go beyond the strictly positive curvature assumptions on $\Om$ and $\dOm$ by circumventing the Reilly formula
via an integration by parts against a carefully chosen test function to obtain  explicit estimates for the crucial bulk-boundary interaction terms also under negative lower curvature bound assumptions.
The second main extension is concerned with an analogous approach to the derivation of an explicit logarithmic Sobolev inequality. We develop the interpolation for this case and show how to obtain the needed boundary-bulk interaction estimates, again via integration by parts against properly chosen test functions. As side results we obtain new explicit geometric estimates of independent interest for the first Steklov eigenvalue (Corollary \ref{cor:stekbound}), the boundary-trace operator for Sobolev functions (Proposition \ref{prop:bdrytraceopest}) and the boundary trace logarithmic Sobolev constant (Proposition \ref{prop:bdrlogsob}), which in special cases was considered before in e.g.\  \cite{beckner,zbMATH05077944}. As we go along we also discuss how to treat the case $\beta =0$ of purely sticky reflection without boundary diffusion.
\bigskip

\section{Poincar\'{e}
Inequality}\label{sec:pi}

For simplicity throughout this paper we  restrict the analysis to the case $\beta=1$ (and $\beta=0$ for the case of pure sticky reflection) and $\gamma>0$, where we assume  that $\Om$ and $\dOm$ have finite (Hausdorff) measure  $|\Om|<\infty $ and $|\dOm|<\infty$. Furthermore by $\lo$ resp. $\ld$ we denote the normalised volume measure on $\Om$ resp. normalised Hausdorff measure on $\dOm$ and choose $\al\in(0,1)$, such that
\begin{equation*}
	\frac{\alpha}{1-\alpha} \frac{|\dOm|}{|\Om|} = \ga.
\end{equation*}
Moreover we set
\begin{equation*}
	\laa := \al \lo + (1-\al)\ld.
\end{equation*}
and find that $-A$ is $\la_\al$-symmetric.\\

Our aim is to estimate the Poincar\'{e}  -- and further down below --  the  logarithmic Sobolev constants for such processes in terms of the geometry of $\Omega$ and its boundary $\dOm$.\\

\noindent We say that a Poincar\'{e} Inequality is fulfilled if there is a constant $C_\al$ such that $\forall f\in C^1(\Om)$
\begin{equation*}
Var_{\laa}(f) \le C_\al \eE_\al(f),	
\end{equation*}
where for $f\in C^1(\Om)$
\begin{gather*}
	Var_{\laa}(f) = \int_\Om f^2 d\laa - \lp\int_\Om f d\laa \rp^2 \\
	\eE_{\al} (f) = \al \int_\Om |\na f|^2 d\lo + (1-\al) \int_{\dOm} |\na^\tau f|^2 d\ld
\end{gather*}
and $\na^\tau$ denotes the tangential derivative operator on $\dOm$.
In the following we denote by $C_\al$ the optimal such constant. By $C_\Om$ and $C_{\dOm}$ we denote the usual (Neumann) Poincar\'{e} constants of $\Om$ and $\dOm$ respectively. We assume that $C_\Om$ and $C_{\dOm}$ (or respective upper bounds for them) are known. In~\cite[Proposition 2.1]{mvr} the following statement was proved in the setting introduced above using an interpolation approach:
\begin{proposition}\label{prop:mvr}
Assume there exist constants $K_{\dOm,\Om},K_1,K_2$ such that for any $f\in C^1(\Om)$
\begin{equation}\label{eq:varboundary}
Var_{\ld} (f) \le K_{\dOm, \Om} \int_{\Om} |\nabla f|^2 d\lo
\end{equation}
and
\begin{equation}\label{eq:mixedterm}
\lp\int_{\Om} f d\lo - \int_{\dOm} f d\ld\rp^2 \le K_1 \int_\Om |\nabla f|^2 d\lo + K_2\int_{\dOm} |\nabla^\tau f|^2 d\ld,
\end{equation}
then it holds for any $\al\in(0,1)$
\begin{align*}
C_\al \le \max\Bigl( &C_\Om + (1-\al)K_1, \al K_2,  \\ & \frac{(1-\al)K_{\dOm,\Omega}C_{\dOm} + \al C_\Om C_{\dOm} + \al(1-\al)(K_{\dOm,\Omega}K_2 + C_{\dOm} K_1)}{(1-\al)K_{\dOm,\Om} + \al C_{\dOm}}\Bigr).
\end{align*}
\end{proposition}
In~\cite[section 3.2]{mvr} constants $K_{\dOm,\Om},K_1,K_2$ were found under the assumption of a positive lower bound on Ricci curvature and a positive lower bound on the second fundamental form on the boundary $\dOm$ (i.e.\ a convex boundary). Our aim is to find $K_{\dOm,\Om},K_1,K_2$ and thus an upper bound on $C_\al$ assuming any upper bound on Sectional curvature and lower bound on Ricci curvature and any upper and lower bound for the second fundamental form on the boundary and to thereby generalise section 3.2 in~\cite{mvr}.\\
The constants $K_1$ and $K_2$ in~\eqref{eq:mixedterm} are not unique. In fact we could add another optimisation over the pairs of values $K_1, K_2$ fulfilling~\eqref{eq:mixedterm} and thus improve the upper bound on $C_\al$.

\smallskip
In order to find suitable constants $K_1, K_2$ fulfilling \eqref{eq:mixedterm} we start with the following proposition where we do not yet make any assumptions on the geometry of the manifold. Later it will be combined with assumptions on curvature and second fundamental form  to obtain explicit admissible choices for $K_1, K_2$.
\begin{proposition}\label{prop:k1k2}
For any $\varphi\in C^1(\Om)$ such that $\frac{\partial \varphi}{\partial N}\vert_{\dOm}=\pm1$ and $\nabla \varphi$ is Lipschitz continuous on $\Om$ inequality~\eqref{eq:mixedterm} in Proposition \ref{prop:mvr} is fulfilled with $K_2=0$ and
\begin{equation*}
K_1 = \lp\frac{|\Om|}{|\dOm|}\rp^2  \inf_{\eps \in (0,\infty)} \left[ (1+\eps)|\nabla \varphi|_2^2 +(1+\eps^{-1})C_\Om |\Delta \varphi|_2^2 \right],
\end{equation*}
where $|\cdot|_2$ denotes the $L^2$-norm on $\Om$ with respect to $\lo$.
\end{proposition}

\begin{proof}
Let $f\in C^1(\Om)$. Without loss of generality we can assume that $\int_{\Om} f d\lo =0$. Now $\forall \eps>0$
\begin{align*}
&\lp \int_\Om f d\lo - \int_{\dOm} f d\ld\rp^2 \\
= &\lp \int_{\dOm} f d\ld\rp^2 = \lp \int_{\dOm} f \frac{\partial \varphi}{\partial N}d\ld \rp^2\\
= &\lp \frac{|\Om|}{|\dOm|} \int_\Om \nabla f \cdot \nabla \varphi d\lo + \frac{|\Om|}{|\dOm|}\int_\Om f\Delta \varphi d\lo\rp^2 \\
\le & \lp \frac{|\Om|}{|\dOm|}\rp^2  \left[ (1+\eps)\lp \int_\Om \nabla f\cdot \nabla\varphi d\lo \rp^2 + (1+\eps^{-1})\lp \int_\Om f\Delta \varphi d\lo \rp^2 \right] \\
\le & \lp \frac{|\Om|}{|\dOm|}\rp^2 \left[ (1+\eps)  \int_\Om |\nabla f|^2 d\lo \int_\Om |\nabla\varphi|^2 d\lo + (1+\eps^{-1}) \int_\Om f^2 d\lo \int_\Om 					(\Delta \varphi)^2 d\lo \right] \\
\le & \lp \frac{|\Om|}{|\dOm|}\rp^2 \left[ (1+\eps)  \int_\Om |\nabla f|^2 d\lo \int_\Om |\nabla\varphi|^2 d\lo + (1+\eps^{-1}) C_\Om \int_\Om |\nabla f|^2
	d\lo \int_\Om (\Delta \varphi)^2 d\lo\right]\\
=& \lp \frac{|\Om|}{|\dOm|}\rp^2 \left[ (1+\eps) \int_\Om |\nabla\varphi|^2 d\lo + (1+\eps^{-1}) C_\Om \int_\Om (\Delta \varphi)^2 d\lo\right]\int_\Om |\nabla f|^2 d\lo.
\end{align*}
\end{proof}

We next find $K_{\dOm,\Om}$ such that inequality~\eqref{eq:varboundary} is fulfilled by proceeding similarly as in the proof of Proposition \ref{prop:k1k2}:

\begin{proposition}\label{prop:kdomom}
For any $\rho\in C^1(\Om)$ such that $\frac{\partial \rho}{\partial N}\vert_{\dOm}=-1$ and $\nabla \rho$ is Lipschitz continuous on $\Om$ inequality~\eqref{eq:varboundary} in Proposition \ref{prop:mvr} is fulfilled with
\begin{equation*}
K_{\dOm,\Om} = \frac{|\Om|}{|\dOm|} \lp 2 |\nabla \rho|_\infty C_\Om^{1/2} + |(\Delta \rho)^-|_\infty C_\Om \rp,
\end{equation*}
where $(\cdot)^-$ denotes the negative part of a function and $|\cdot|_\infty$ denotes the $L^\infty$-norm on $\Om$ with respect to $\lo$.
\end{proposition}

\begin{proof}
Let $f\in C^1(\Om)$. Without loss of generality we can assume that $\int_{\Om} f d\lo =0$. Now we can calculate similarly as in the previous result
\begin{align*}
Var_{\ld}(f)&=\int_{\dOm} f^2 d\ld -\lp \int_{\dOm} f d\ld \rp^2 \le \int_{\dOm} f^2 d\ld = - \int_{\dOm} f^2 \frac{\partial \rho}{\partial N} d\ld \\
	&= - \frac{|\Om|}{|\dOm|} \int_\Om 2f \nabla f\cdot \nabla \rho d\lo - \frac{|\Om|}{|\dOm|} \int_\Om f^2\Delta \rho d\lo\\
	&\le 2 \frac{|\Om|}{|\dOm|}  \int_\Om |f||\nabla f||\nabla \rho| d\lo + \frac{|\Om|}{|\dOm|}  \int_\Om f^2 (\Delta \rho)^- d\lo\\
	&\le 2 \frac{|\Om|}{|\dOm|}  |\nabla \rho|_\infty \int_\Om |f||\nabla f| d\lo + \frac{|\Om|}{|\dOm|}  |(\Delta \rho)^-|_\infty \int_\Om f^2 d\lo\\
	&\le 2 \frac{|\Om|}{|\dOm|}  |\nabla \rho|_\infty \lp \int_\Om f^2 d\lo \int_\Om |\nabla f|^2 d\lo\rp^{1/2} + \frac{|\Om|}{|\dOm|}  |(\Delta \rho)^-|_\infty \int_\Om f^2 d\lo \\
	&\le 2 \frac{|\Om|}{|\dOm|}  |\nabla \rho|_\infty C_\Om^{1/2} \int_\Om |\nabla f|^2 d\lo + \frac{|\Om|}{|\dOm|}  |(\Delta \rho)^-|_\infty C_\Om \int_\Om |\nabla f|^2 d\lo \\
	&=\frac{|\Om|}{|\dOm|} \lp 2 |\nabla \rho|_\infty C_\Om^{1/2} + |(\Delta \rho)^-|_\infty C_\Om \rp \int_\Om |\nabla f|^2 d\lo.
\end{align*}
\end{proof}

\begin{remark}\label{rem:steklov}
Denote by $\sigma$ the first nontrivial eigenvalue of the Steklov eigenvalue problem
\begin{equation*}
	\begin{cases}\Delta f = 0, &\text{ in } \Om\\
		\frac{\pd f}{\pd N} = \sigma f, &\text{ on } \dOm,
	\end{cases}
\end{equation*}
which is characterised (using normalised measures) by
\begin{equation*}
\sigma = \frac{|\Om|}{|\dOm|} \inf_{\substack{f\in C^1(\Om)\\ \int_{\dOm} f d\ld = 0}}  \frac{\int_\Om |\nabla f|^2 d\lo}{\int_{\dOm} f^2 d\ld}.
\end{equation*}
Thus we have for the optimal constant $K_{\dOm,\Om}$ in inequality~\eqref{eq:varboundary} in Proposition \ref{prop:mvr}
\begin{equation*}
K_{\dOm,\Om} = \sup_{f\in C^1(\Om)} \frac{Var_{\ld}(f)}{\int_\Om |\nabla f|^2 d\lo} = \sup_{\substack{f\in C^1(\Om)\\ \int_{\dOm} f d\ld =0}} \frac{\int_{\dOm} f^2 d\ld}{\int_\Om |\nabla f|^2 d\lo}
= \frac{|\Om|}{|\dOm|} \sigma^{-1}.
\end{equation*}
Therefore by finding upper bounds for the optimal $K_{\dOm,\Om}$ we find lower bounds for the first nontrivial Steklov eigenvalue. We use this connection in the computations for Example \ref{ex:hyperbolic} below.
\end{remark}

To obtain an explicit constant it now remains to specify functions $\varphi$ and $\rho$ with the desired properties. Note that despite fulfilling the same assumptions, $\varphi$ and $\rho$ may be chosen independently in order to optimise the estimates. It seems natural to define both functions of the form $\psi \circ \rho_{\dOm}$ for some appropriate function $\psi$, where $\rho_{\dOm}$ denotes the distance to the boundary function.\\
We use for $k,\gamma\in\mathbb{R}$ the function
\begin{equation}\label{eq:hfctn}
h:[0,\infty)\to\mathbb{R}, \ h(t):=\begin{cases}
	cos(\sqrt{k}t)-\frac{\gamma}{\sqrt{k}}sin(\sqrt{k}t), & k\ge0 \\
	cosh(\sqrt{-k}t)-\frac{\gamma}{\sqrt{-k}}sinh(\sqrt{-k}t), & k<0.
\end{cases}
\end{equation}
Let $h^{-1}(0):=\inf\{t\ge 0: h(t)=0\}$, where $h^{-1}(0)=\infty$ if $h(t)>0$ for all $t\ge 0$.
We denote by $Ric$ and $sect$ the Ricci and sectional curvatures of $\Om$,  and by $\Pi$ the second fundamental form on the boundary $\dOm$, i.e.
\begin{equation*}
\Pi(X,Y):=\langle \nabla_X N, Y \rangle, \ X,Y\in T_x \partial \Om, x\in\partial \Om,
\end{equation*}
where $N$ is the outward pointing unit normal vector field of $\partial \Om$.

\begin{lemma}\label{lem:rhol2}
Let $k_1,k_2\in\mathbb{R}$ such that $Ric\ge k_1 (d-1),  sect \le k_2$ and $\gamma_1,\gamma_2\in\mathbb{R}$ such that $\gamma_1 id \le \Pi \le \gamma_2 id$. We construct a function  $\varphi\in C^1(\Om)$ such that $\frac{\partial \varphi}{\partial N}\vert_{\dOm}=-1$ and $\nabla \varphi$ is Lipschitz continuous on $\Om$ and for $t_0\in (0,h_2^{-1}(0))$
\begin{align*}
|\nabla \varphi|_2^2&\le \frac{1}{|\Om|} \int_0^{t_0} H_{d-1}(\{\rho_{\dOm}=t\}) \lp 1-\frac{t}{t_0}\rp^2 dt, \\
|\Delta \varphi|_2^2 &\le \frac{1}{|\Om|} \int_0^{t_0} H_{d-1}(\{\rho_{\dOm} = t\}) \biggl( \lp \lp (d-1)\frac{h_2'}{h_2}(t)  \lp 1-\frac{t}{t_0} \rp - \frac{1}{t_0}\rp^- \rp^2  \\
	&+ \lp \lp (d-1)\frac{h_1'}{h_1}(t)  \lp 1-\frac{t}{t_0} \rp -\frac{1}{t_0}\rp^+\rp^2 \biggr) dt,
\end{align*}
where $h_i, i=1,2$ are as defined above in~\eqref{eq:hfctn} with $k=k_i$ and $\gamma=\gamma_i$.
\end{lemma}

\begin{proof} It is easy to see that $h_2^{-1}(0)\le h_1^{-1}(0).$   Let $\rho_{\dOm}$ be the distance function to the boundary. By the Laplacian comparison theorem, we have
\begin{align}
\Delta \rho_{\dOm} \le \frac{(d-1)h_1'}{h_1}(\rho_{\dOm})\text { on }  \{\rho_{\dOm} < h_1^{-1}(0)\}, \label{eq:laplacecomp1}\\
\Delta \rho_{\dOm} \ge \frac{(d-1)h_2'}{h_2}(\rho_{\dOm})\text{ on } \{\rho_{\dOm} < h_2^{-1}(0)\} \label{eq:laplacecomp2}.
\end{align}
Indeed, \eqref{eq:laplacecomp2} follows from \cite[Theorem 3.1]{wang} for the Laplacian comparison theorem due to \cite{kasue}, and by \cite[Corollary 3.2]{wang}  which says that the injectivity radius of $\dOm$ is larger than $  h_2^{-1}(0)$.
Next, for $x\in \Omega$ with $\rho_{\dOm} (x)<  h_1^{-1}(0)$, let $p\in \partial\Omega$ be the projection such that
$\gamma(s):= \exp[-s N(p)], s\in [0,\rho_{\dOm}(x)]$ be the minimal geodesic from $p$ to $x$. Let
$\{X_i\}_{1\le i\le d-1}$ be orthonormal vector fields around $x$  orthogonal to $\nabla \rho_{\dOm}(x).$
Let $J_i(s))_{s\in [0,\rho_{\dOm}(x)]}$ be the Jacobi fields along  the geodesic $\gamma$ such that
$J_i(\rho_{\dOm}(x))=X_i(x)$ and
$$\langle\dot J_i(0),v\rangle=-\Pi(J_i(0), v),\ \ v\in T_p\partial \Omega.$$
Let $\mathcal R$ be the Riemannian curvature tensor. By the second variational formula (see page 321 in \cite{chavel95}) we have
$${\rm Hess}_{\rho_{\dOm}}(X_i,X_i)(x)=-\Pi(J_i(0), J_i(0))+ \int_0^{\rho_{\dOm}(x)}
\Big(|\dot J_i(s)|^2 -\big\langle\mathcal R(\dot \gamma(s), J_i(s))\dot\gamma, J_i(s)\big\rangle \Big) d s.$$
Let $(X_i(s))_{s\in [0, \rho_{\dOm}(x)]}$ be the parallel displacement of $(X_i(0):=X_i(x))_{s\in [0, \rho_{\dOm}(x)]},$ and denote
$$\tilde J_i(s) = \frac {h_1(s)}{h_1(\rho_{\dOm}(x))}X_i(s),\ \ 1\le i\le d-1, s \in [0, \rho_{\dOm}(x)].$$
 Then the index lemma   yields
$${\rm Hess}_{\rho_{\dOm}}(X_i,X_i)(x)\le -  \Pi(\tilde  J_i(0), \tilde J_i(0))+ \int\limits_0^{\rho_{\dOm}(x)}
\big(|\dot{ \tilde {J}}_i(s)|^2 -\big\langle\mathcal R(\dot \gamma(s), \tilde J_i(s))\dot\gamma, \tilde J_i(s)\big\rangle \big) d s.$$ Noting that $h_1''(s)=-k_1 h_1(s)$, this implies
\eqref{eq:laplacecomp1}. \\
Now for $t_0\in (0,h_2^{-1}(0))$ we define
\begin{equation*}
\varphi= \int_0^{\rho_{\dOm}} \lp 1-\frac{s}{t_0}\rp^+ ds.
\end{equation*}
We have
\begin{equation*}
\nabla \varphi (x)= \begin{cases} \nabla \rho_{\dOm}(x)\cdot \lp 1-\frac{\rho_{\dOm}(x)}{t_0} \rp,  &\rho_{\dOm}(x) \le t_0 \\
	0,  &\text{ else},
	\end{cases}
\end{equation*}
and thus $\frac{\partial \varphi}{\partial N}\vert_{\dOm}=-1$ and $\nabla \varphi$ is Lipschitz continuous. Furthermore
\begin{equation*}
\Delta \varphi (x)= \begin{cases}\Delta \rho_{\dOm}(x) \lp 1- \frac{\rho_{\dOm}(x)}{t_0}\rp -\frac{1}{t_0}, & \rho_{\dOm}(x)\le t_0, \\
	0 &\text{ else.}
	\end{cases}
\end{equation*}
Using the Coarea formula we get
\begin{align*}
\int_\Om |\nabla \varphi|^2 d\lo &= \int_{\{\rho_{\dOm}\le t_0\}} \lp 1-\frac{\rho_{\dOm}}{t_0}\rp^2 |\nabla \rho_{\dOm}|^2 d\lo =  \int_{\{\rho_{\dOm}\le t_0\}} \lp 1-\frac{\rho_{\dOm}}{t_0}\rp^2 d\lo \\
	&=\frac{1}{|\Om|} \int_0^{t_0} \int_{\{\rho_{\dOm}=t\}} \lp 1-\frac{t}{t_0}\rp^2 dH_{d-1} dt \\
	&= \frac{1}{|\Om|} \int_0^{t_0} H_{d-1}(\{\rho_{\dOm}=t\}) \lp 1-\frac{t}{t_0}\rp^2 dt,
\end{align*}
where $H_{d-1}$ denotes the $(d-1)$-dimensional Hausdorff measure. Furthermore
\begin{equation*}
\int_\Om (\Delta \varphi)^2 d\lo = \int_\Om \lp (\Delta \varphi)^+ - (\Delta \varphi)^-  \rp^2 d\lo = \int_\Om \lp (\Delta \varphi)^+\rp^2 d\lo + \int_\Om \lp (\Delta \varphi)^-\rp^2 d\lo.
\end{equation*}
We see that by inequalities~\eqref{eq:laplacecomp1} and~\eqref{eq:laplacecomp2} on $\{\rho_{\dOm}\le t_0\}$
\begin{align*}
\Delta\varphi \ge  (d-1)\frac{h_2'}{h_2}(\rho_{\dOm}) \lp 1-\frac{\rho_{\dOm}}{t_0} \rp -\frac{1}{t_0} &\Rightarrow (\Delta\varphi)^- \le  \lp  (d-1)\frac{h_2'}{h_2}(\rho_{\dOm}) \lp 1-\frac{\rho_{\dOm}}{t_0} \rp -\frac{1}{t_0}\rp^-,\\
\Delta\varphi \le  (d-1)\frac{h_1'}{h_1}(\rho_{\dOm}) \lp 1-\frac{\rho_{\dOm}}{t_0} \rp -\frac{1}{t_0} &\Rightarrow (\Delta\varphi)^+ \le \lp (d-1)\frac{h_1'}{h_1}(\rho_{\dOm})   \lp 1-\frac{\rho_{\dOm}}{t_0} \rp -\frac{1}{t_0}\rp^+.
\end{align*}
Thus
\begin{align*}
\int_\Om (\Delta \varphi)^2 d\lo \le &\int_{\{\rho_{\dOm}\le t_0\}} \lp \lp (d-1)\frac{h_2'}{h_2}(\rho_{\dOm}) \lp 1-\frac{\rho_{\dOm}}{t_0} \rp - \frac{1}{t_0}\rp^- \rp^2 \\
		&+ \lp \lp (d-1)\frac{h_1'}{h_1}(\rho_{\dOm})  \lp 1-\frac{\rho_{\dOm}}{t_0} \rp -\frac{1}{t_0}\rp^+\rp^2 d\lo \\
	=& \frac{1}{|\Om|} \int_0^{t_0} H_{d-1}(\{\rho_{\dOm} = t\}) \biggl( \lp \lp (d-1)\frac{h_2'}{h_2}(t)  \lp 1-\frac{t}{t_0} \rp - \frac{1}{t_0}\rp^- \rp^2  \\
	&+ \lp \lp (d-1)\frac{h_1'}{h_1}(t)  \lp 1-\frac{t}{t_0} \rp -\frac{1}{t_0}\rp^+\rp^2 \biggr) dt.
\end{align*}
\end{proof}

In the previous Lemma $t_0\in (0,h_2^{-1}(0))$ may be chosen to either optimise $|\nabla \varphi|_2^2$ or $|\Delta \varphi|_2^2$.

\begin{lemma}\label{lem:rho}
Let $k_2\in\mathbb{R}$ such that $sect \le k_2$ and $\gamma_2\in\mathbb{R}$ such that $\Pi\le \gamma_2 id$. Then $k_2>-\gamma_2^2$, and for all $\eps>0$ there exists a function $\rho\in C^1(\Om)$ such that $\frac{\partial \rho}{\partial N}\vert_{\dOm}=-1$ and $\nabla \rho$ is Lipschitz continuous on $\Om$ and
\begin{align*}
|\nabla \rho|_\infty &\le 1, \\
|(\Delta \rho)^-|_\infty &\le   \inf_{t_1\in(0,h_2^{-1}(0))} \sup_{t\in(0,t_1)} \lp(d-1)\frac{h_2'}{h_2}\lp t\rp \lp 1-\frac{t}{t_1}\rp - \frac{1}{t_1}\rp^- +\eps.
	 \end{align*}
\end{lemma}

\begin{proof}
Let  $h_2$ be the function defined as in~\eqref{eq:hfctn} with $k=k_2$ and $\gamma=\gamma_2$.  If  $k_2\le -\gamma_2^2,$ then $h_2^{-1}(0)=\infty$, so that  \cite{wang}[Corollary 3.2] implies that the cut locus of $\partial \Omega$ is empty, which is contradictive to the fact that the maximum point of $\rho_{\dOm}$ is in the cut locus.  Hence,   $k_2>-\gamma_2^2.$

\smallskip
Let $t_1\in(0,h_2^{-1}(0))$ to be chosen later.
By~\eqref{eq:laplacecomp2}
\begin{equation}\label{eq:laplacea}
\Delta \rho_{\dOm} \ge (d-1)\frac{h_2'}{h_2}(\rho_{\dOm})
 \text{ on } \{\rho_{\dOm} \le t_1\}.
\end{equation}
Now define
\begin{equation*}
\rho= \int_0^{\rho_{\dOm}} \lp 1-\frac{s}{t_1}\rp^+ ds.
\end{equation*}
We have
\begin{equation*}
\nabla \rho (x)= \begin{cases} \nabla \rho_{\dOm}(x)\cdot \lp 1-\frac{\rho_{\dOm}(x)}{t_1} \rp,  &\rho_{\dOm}(x) \le t_1 \\
	0,  &\text{ else},
	\end{cases}
\end{equation*}
and thus $\frac{\partial \rho}{\partial N}\vert_{\dOm}=-1$, $|\nabla \rho|_\infty\le 1$ and $\nabla \rho$ is Lipschitz continuous. Furthermore
\begin{equation*}
\Delta \rho (x)= \begin{cases}\Delta \rho_{\dOm}(x) \lp 1- \frac{\rho_{\dOm}(x)}{t_1}\rp -\frac{1}{t_1}, & \rho_{\dOm}(x)\le t_1\\
	0, & \text{ else},
	\end{cases}
\end{equation*}
and thus by inequality~\eqref{eq:laplacea} on $\{\rho_{\dOm}\le t_1\}$
\begin{equation*}
\Delta \rho \ge (d-1)\frac{h_2'}{h_2}\lp \rho_{\dOm}\rp \lp 1 - \frac{\rho_{\dOm}}{t_1}\rp - \frac{1}{t_1} \Rightarrow (\Delta \rho)^- \le \lp(d-1)\frac{h_2'}{h_2}\lp \rho_{\dOm}\rp \lp 1 - \frac{\rho_{\dOm}}{t_1}\rp - \frac{1}{t_1}\rp^-.
\end{equation*}
We can still choose $t_1\in(0,h_2^{-1}(0))$ to obtain for arbitrary $\eps>0$
\begin{equation*}
|(\Delta \rho)^-|_\infty \le \inf_{t_1\in(0,h_2^{-1}(0))} \sup_{t\in(0,t_1)} \lp(d-1)\frac{h_2'}{h_2}\lp t\rp \lp 1-\frac{t}{t_1}\rp - \frac{1}{t_1}\rp^- + \eps.
\end{equation*}
 \end{proof}

Inserting $\varphi$ and $\rho$ as defined in Lemma \ref{lem:rhol2} and Lemma \ref{lem:rho} in Proposition \ref{prop:k1k2} and Proposition \ref{prop:kdomom} we now get explicit constants $K_1, K_2$ and $K_{\dOm,\Om}$ in terms of bounds on sectional and Ricci curvature and second fundamental form on the boundary. We state these in the following Proposition.

\begin{proposition}\label{prop:explconst}
Let $k_1,k_2\in\mathbb{R}$ such that $Ric \ge (d-1)k_1,  sect \le k_2$ and $\gamma_1,\gamma_2\in\mathbb{R}$ such that $\gamma_1 id \le \Pi \le \gamma_2 id$. Then the assumptions in Proposition \ref{prop:mvr} are fulfilled with
\begin{align*}
K_1 = \frac{|\Om|}{|\dOm|^2}  \inf_{t_0\in(0,h_2^{-1}(0))} \inf_{\eps \in (0,\infty)} &\biggl[\int_0^{t_0}(1+\eps) H_{d-1}(\{\rho_{\dOm}=t\}) \lp 1-\frac{t}{t_0}\rp^2 \\ &
+\lp1+\frac{1}{\eps}\rp C_\Om H_{d-1}(\{\rho_{\dOm} = t\}) \biggl( \lp \lp (d-1)\frac{h_2'}{h_2}(t)  \lp 1-\frac{t}{t_0} \rp - \frac{1}{t_0}\rp^- \rp^2 \\
&+ \lp \lp (d-1)\frac{h_1'}{h_1}(t)  \lp 1-\frac{t}{t_0} \rp -\frac{1}{t_0}\rp^+\rp^2 \biggr) dt \biggr],
\end{align*}
\begin{gather*}
K_2 = 0,\\
K_{\dOm,\Om} =    \frac{|\Om|}{|\dOm|} \lp 2 C_\Om^{1/2} + C_\Om \inf_{t_1\in(0,h_2^{-1}(0))} \sup_{t\in(0,t_1)} \lp(d-1)\frac{h_2'}{h_2}\lp t\rp \lp 1-\frac{t}{t_1}\rp - \frac{1}{t_1}\rp^-\rp.
\end{gather*}
\end{proposition}

As explained in Remark \ref{rem:steklov}, upper bounds for the optimal $K_{\dOm,\Om}$ correspond to lower bounds for the first nontrivial Steklov eigenvalue. Thus we now also get a lower bound on the first nontrivial Steklov eigenvalue $\sigma$ that is explicit in terms of upper bounds on sectional curvature and second fundamental form on the boundary:

\begin{corollary} \label{cor:stekbound}
Let $k_2\in\mathbb{R}$ such that $sect \le k_2$ and $\gamma_2\in\mathbb{R}$ such that $\Pi\le \gamma_2 id$. Then for the first non-trivial Steklov eigenalue $\sigma$ of $\Om$ it holds that
\begin{equation*}
\sigma \ge    \lp 2 C_\Om^{1/2} + C_\Om \inf_{t_1\in(0,h_2^{-1}(0))} \sup_{t\in(0,t_1)} \lp(d-1)\frac{h_2'}{h_2}\lp t\rp \lp 1-\frac{t}{t_1}\rp - 	\frac{1}{t_1}\rp^-\rp^{-1}.
\end{equation*}
\end{corollary}

\begin{example}
As an example one may consider a ball of radius one around the origin in the hyperbolic plane. Then $C_\Om\approx 0.3377$ and $h_2(t)=cosh(t)-coth(1)sinh(t)$ (for more details on this see also Example \ref{ex:hyperbolic} below). This results in the lower bound $\sigma \ge 0.5145$ while as stated in~\cite{escobar} $\sigma = coth(1)-tanh(1/2) \approx 0.8509$ for this example.
\end{example}

By inserting the set of constants stated in Proposition \ref{prop:explconst} into Proposition \ref{prop:mvr} we get an explicit upper bound on the Poincar\'{e} constant again in terms of bounds on sectional and Ricci curvature and second fundamental form. For comparison we state the following obvious upper bound for the optimal Poincar\'{e} constant without the interpolation, which can be derived from Proposition \ref{prop:mvr} by considering the limit as $K_{\dOm,\Om}$ tends to infinity.

\begin{proposition}\label{prop:trivialpc}
Assume there exist constants $K_1,K_2$ such that for any $f\in C^1(\Om)$
	\begin{equation*}
	\lp\int_{\Om} f d\lo - \int_{\dOm} f d\ld\rp^2 \le K_1 \int_\Om |\nabla f|^2 d\lo + K_2\int_{\dOm} |\nabla^\tau f|^2 d\ld,
	\end{equation*}
	then it holds for any $\al\in(0,1)$
	\begin{equation*}
		C_\al \le \max \lp C_\Om + (1-\al)K_1, C_{\dOm} + \al K_2\rp.
	\end{equation*}
\end{proposition}

\subsection{Examples}

In the following we consider as examples balls in Euclidean plane and in hyperbolic plane and compare the results obtained above with or without the interpolation approach. While the former example has already been treated in~\cite{mvr}, the latter was not included there due to negative curvature.

\subsubsection{Ball in $
\mathbb R^2$} \label{ex:euclidean}

Let $\Om := B_1$ be the unit ball in $\mathbb{R}^2$. In this case the sectional curvature equals $k=0$, and for the second fundamental form on the boundary we have $\gamma=1$. The constants $K_1, K_2$ and $K_{\dOm,\Om}$ are now to be computed for this specific example. The general results in~\cite[section 3.2]{mvr} concerning the values for these constants are not applicable. However the same paper also contains a computation adapted to this specific example. In the following we will consider the values for the three constants obtained by computations adapted to this specific example as well as the values obtained by the above general results. For both sets of constants we will compare the upper bounds on Poincar\'{e} constants obtained by the interpolation approach with the bounds obtained without interpolation as in Proposition \ref{prop:trivialpc} and the exact values for the Poincar\'{e} constants.

\begin{figure}
	\includegraphics[width=\linewidth]{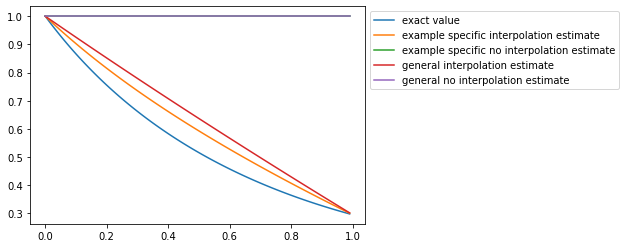}
	\caption{Exact Poincar\'{e} constants (blue), interpolation (yellow) and no interpolation (green) results using computations specifically adapted to the example, interpolation (red) and no interpolation (purple) results using computations not specifically adapted to the example.}
	\label{fig:r2ball}
\end{figure}

\smallskip
We first recall the results of the computations adapted to the ball example made in~\cite[section 3.1]{mvr}:
\begin{equation*}
C_\Om \approx \frac{1}{3.39},\ C_{\dOm}=1,\ K_1=\frac{3}{16},\ K_2=0,\ K_{\dOm,\Om}=\frac{1}{2}.
\end{equation*}
The upper bound obtained from Proposition \ref{prop:mvr} is
\begin{equation*}
C_\al \le \frac{8(1-\al) + 16\al C_\Om + 3\al(1-\al)}{8(1+\al)},
\end{equation*}
while the upper bound obtained for the same values of $K_1,K_2$ and $K_{\dOm,\Om}$ from Proposition \ref{prop:trivialpc} is $C_\al \le 1$. Moreover we also refer to~\cite{mvr} for the procedure to calculate the exact values for $C_\al, \al\in(0,1).$ Using the results from the previous pages instead, we find different constants: We have $h_1(t)=h_2(t)=1-t$ from which follows by Lemma \ref{lem:rhol2} and Lemma \ref{lem:rho} that
\begin{equation*}
\forall \eps >0\ \exists \varphi: |\nabla \varphi|_2^2 \le \eps,\ |\Delta \varphi|_2^2\le1 + \eps
\end{equation*}
and
\begin{equation*}
\forall \eps >0\ \exists \rho:\ |\nabla \rho|_{\infty}\le 1,\ |(\Delta \rho)^-|_{\infty}\le 2 + \eps.
\end{equation*}
Inserting this in Proposition \ref{prop:k1k2} and Proposition \ref{prop:kdomom} we get
\begin{equation*}
K'_1 = \frac{C_\Om}{4} ,\ K'_2=0,\ K'_{\dOm,\Om}\approx 0.8381.
\end{equation*}
Inserting this in Proposition \ref{prop:trivialpc} results in $C_\al \le 1$ while inserting in Proposition \ref{prop:mvr} we get
\begin{equation*}
C_\al \le \max \lp C_\Om + (1-\al)K'_1, \frac{(1-\al)K'_{\dOm,\Om}+\al C_\Om + \al(1-\al)K'_1}{(1-\al)K'_{\dOm,\Om}+\al}\rp.
\end{equation*}
We depict these results in Figure \ref{fig:r2ball}. Note that the green and purple curves overlap.
From this we see that the upper bounds obtained from the above general results are only slightly worse than the upper bounds obtained by computing with the specific example in mind. Furthermore it is obvious from the proof of Proposition \ref{prop:mvr}, that results obtained from the interpolation approach must be at least as good as results obtained without interpolation. However the figure shows
for both sets of constants that the interpolation results clearly differ from the no interpolation results and give significantly better bounds than the approach without interpolation.\\
We remark that using the interpolation approach with $K_1', K_2'$ and $K_{\dOm,\Om}$ (i.e. combining the respective best constants from above) would result in an even better approximation of the curve of exact values.

\subsubsection{Ball in hyperbolic plane} \label{ex:hyperbolic}
We consider the unit metric ball in the hyperbolic plane and compute the exact Poincar\'{e} constants $C_\alpha, \alpha\in(0,1)$ numerically. As in the previous example we compare these with general/specific example and interpolation/no interpoaltion results. In more detail we consider the unit ball $B_1(0)\subset \mathbb{R}^2$ with the hyperbolic metric
\begin{equation}\label{eq:hypmet}
	g_h = \frac{4}{(1-|x|^2)^2}g,
\end{equation}
where $g=(dx^1)^2+ (dx^2)^2$ is the standard metric in $\mathbb{R}^2$, resulting in the space form of constant sectional curvature $K=-1$. $\Om$ will be a unit ball in this hyperbolic plane. We will start by computing the exact values for $C_\al, \al\in(0,1).$
Using that $|\Om|=2\pi (cosh(1)-1)$ and $|\dOm|=2\pi sinh(1)$ the operator $A=A_\alpha$ associated with the Dirichlet form $\eE_{\alpha}$ then becomes
\begin{equation*}
	A_\al f = \Del f\1_\Om + \lp \Del^\tau f - \frac{\al}{1-\al} \frac{sinh(1)}{cosh(1)-1} \frac{\pd f}{\pd N}\rp\1_{\dOm}.
\end{equation*}
An eigenvector of $-A_\al$ for eigenvalue $\la\ge0$ is then a function $f\in D(A_\al)$ such that the following system of partial differential equations is fulfilled
\begin{align*}\begin{cases}
		\Del f = -\la f  & \text{ in } \Om\\
		\Del^\tau f - \frac{\al}{1-\al} \frac{sinh(1)}{cosh(1)-1} \frac{\pd f}{\pd N} = - \la f & \text{ on } \dOm
	\end{cases}.
\end{align*}
Since $f$ and $A_\al f$ are continuous on $\Om$, this is equivalent to
\begin{align}\label{eq:hypsyst}
	\begin{cases}
		\Del f = -\la f  & \text{ in } \Om\\
		\Del^\tau f - \frac{\al}{1-\al} \frac{sinh(1)}{cosh(1)-1} \frac{\pd f}{\pd N} =  \Del f & \text{ on } \dOm
	\end{cases}.
\end{align}
Following the well-known procedure for the Laplacian with Neumann boundary conditions, see e.g.~\cite[Chapter 2.5]{chavel},
we introduce spherical coordinates about $x=0$ by
\begin{equation*}
	x = r\xi,\ r=tanh(t/2),
\end{equation*}
where $r\in[0,1],\ t\in[0,\infty),\ \xi\in\mathbb{S}^1$. $\Om$ is then characterised by restriction of $t$ to $[0,1]$. In these coordinates~\eqref{eq:hypmet} reads
\begin{equation*}	g_h = (dt)^2 + sinh^2(t)|d\xi|^2.
\end{equation*}
We then separate variables, i.e. $f(t,\xi)=T(t)G(\xi)$. Furthermore we denote by $\square$ the Laplacian on $\mathbb{S}^1$ and by $'$ differentiation with respect to $t$.
Using the Laplacian in spherical coordinates, the first equation of~\eqref{eq:hypsyst} becomes
\begin{equation*}
	sinh(t)^{-1}(sinh(t)T'(t))'G(\xi) + sinh(t)^{-2}T(t)\square G(\xi) = - \lambda T(t)G(\xi).
\end{equation*}
and in terms of $G$ and $T$
\begin{align*}
	\begin{cases}
		\square G(\xi) + \gamma G(\xi) =0\\
		(sinh(t)T'(t))' + (\lambda-\gamma sinh(t)^{-2})sinh(t)T(t) =0,
	\end{cases}
\end{align*}
where $\gamma= l^2,\ l\in \mathbb{N}$ are the negatives of the eigenvalues of $\square$ on $\mathbb{S}^1$. According to~\cite[Chapter 12.5]{chavel} the solution $T$ is given via
\begin{equation*}
	T(t)= P^\mu_\nu(cosh(t)),
\end{equation*}
where $P^\mu_\nu(\cdot)$ is the associated Legendre function of first kind with $\mu$ and $\nu$ given via
\begin{equation*}
	\mu =l, \ \nu= -\frac{1}{2} \pm \sqrt{-\lambda + \frac{1}{4}}.
\end{equation*}
We thus obtain a two parameter family of eigenfunctions $f_{n,l}(t,\xi)=P^l_n(cosh(t))G_l(\xi)$ where $G_l$ is the eigenfunction for eigenvalue $-l^2$ and $n$ respectively $\lambda_n$ is constrained via the boundary condition as follows. Using that $\Del^\tau f = \frac{1}{sinh^2(1)}T\square G$, the second equation of~\eqref{eq:hypsyst} which holds on the boundary, i.e. for $t=1$, becomes
\begin{align*}
	&\Del (TG)(1,\xi) = \frac{1}{sinh^2(1)} T(1)\square G(\xi) - \frac{sinh(1)}{cosh(1)-1}\frac{\al}{1-\al} T'(1)G(\xi) \\
	\Leftrightarrow &T''(1) + T'(1)\lp\frac{cosh(1)}{sinh(1)} + \frac{sinh(1)}{cosh(1)-1}\frac{\al}{1-\al}\rp = 0 \\
	\Leftrightarrow &{P^l_\nu}'' (cosh(1))sinh^2(1) + {P^l_\nu}' (cosh(1))\lp 2cosh(1) + \frac{sinh^2(1)}{cosh(1)-1}\frac{\al}{1-\al} \rp =0.
\end{align*}
We obtain a two-parameter family of values $\lambda_{l,n}$ satisfying this equation. For $\al\in(0,1)$ $\lambda_\al := \min_{l,n} \lambda_{l,n}$ is the desired spectral gap and $C_\alpha=1/\lambda_\al$.\\

As we found no explicit account of $C_\Om$ and $C_{\dOm}$ we also compute these here. Following the same procedure as above including spherical coordinates and separation of variables we see that eigenfunctions $f$ of the Laplacian on $\Om$ with Neumann boundary conditions on $\dOm$ are again of the form
\begin{equation*}
	f(t,\xi)=T(t)G(\xi), \text{ with } T(t)=P^\mu_\nu (cosh(t)),\ \mu=l,\ \nu=-\frac{1}{2} \pm \sqrt{-\la + \frac{1}{4}},
\end{equation*}
where $G$ are eigenfunctions of the Laplacian on $\mathbb{S}^1$ for the eigenvalues $-l^2,\ l\in\mathbb{N}$.\\
Now the boundary condition amounts to
\begin{equation*}
	\frac{\pd f}{\pd N} = 0 \text{ on } \dOm\  \Leftrightarrow\ {P^\mu_\nu}' (cosh(1)) = 0.
\end{equation*}
We again obtain a two-parameter family of values $\lambda_{l,n}$ satisfying this equation. Again $\lambda_1 := \min_{l,n} \lambda_{l,n}\approx 2.9614$ is the desired spectral gap and $C_\Om=1/\lambda_1\approx 0.3377$.


Furthermore by scaling we can deduce that since the eigenvalues of $\Delta_{S^1_\mathbb{R}}$ are known to be $\lambda_k := -k^2,\ k\in\mathbb{N}_0$ (see e.g.~\cite[Chapter 2.4]{chavel}), those of $\Delta_{S^1_H}$ are $\tilde{\lambda}_k= -\frac{k^2}{sinh^2(1)},\ k\in\mathbb{N}_0$ and thus $C_{\dOm} = sinh^2(1)$.\\

We now compute $K_1$ and $K_2$ in a fashion adapted to the specific example. The following computation is similar to the one referenced above in Example \ref{ex:euclidean}.
Using that for $f\in L^1(\dOm)$:
\begin{equation*}
	\int_{\dOm} f(y) \ld (dy) = \int_{\Om} f(x/|x|) \lo(dx),
\end{equation*}
we get
\begin{align*}
	&\lp \int_\Om f d\lo - \int_{\dOm} f d\ld \rp^2 \\
	&\le \int_{\Om} (f(x)-f(x/\left|x\right|)^2 \lo(dx)\\
	&=  \frac{1}{|\Om|} \int_{\dOm} \int_0^1 \lp f(tanh(t/2)\xi) - f(tanh(1/2)\xi) \rp^2 sinh(t) dt d\xi \\
	&=  \frac{1}{|\Om|} \int_{\dOm} \int_0^1 \lp \int_t^1 \frac{d}{ds} f(tanh(s/2)\xi)ds \rp^2 sinh(t) dt d\xi \\
	&\le  \frac{1}{|\Om|} \int_{\dOm} \int_0^1 (1-t) \int_t^1 \lp \frac{d}{ds} f(tanh(s/2)\xi)\rp^2 ds\ sinh(t) dt d\xi \\
	&=  \frac{1}{|\Om|} \int_{\dOm} \int_0^1 \int_0^s (1-t)sinh(t) dt \lp \frac{d}{ds} f(tanh(s/2)\xi)\rp^2 ds d\xi \\
	&=  \frac{1}{|\Om|} \int_{\dOm} \int_0^1 \lp sinh(s) - (s-1)cosh(s) -1\rp \lp \langle \nabla f (tanh(s/2)\xi),\frac{d}{ds}tanh(s/2)\xi \rangle \rp^2 ds d\xi \\
	&\le  \frac{1}{|\Om|} \int_{\dOm} \int_0^1 \lp sinh(s) - (s-1)cosh(s) -1\rp |\nabla f (tanh(s/2)\xi)|^2 ds d\xi \\
	&\le K_1 \frac{1}{|\Om|} \int_{\dOm} \int_0^1 |\nabla f (tanh(s/2)\xi)|^2 sinh(s) ds d\xi \\
	& = K_1 \int_{\Om} |\nabla f (x)|^2 \lo (dx),
\end{align*}
where
\begin{equation*}
	K_1 := \max_{s\in[0,1]} \frac{sinh(s)-(s-1)cosh(s) -1}{sinh(s)} \approx 0.1782.
\end{equation*}
In particular, with $K_1$ as above setting $K_2=0$ is an admissible choice.

\smallskip
As explained in Remark \ref{rem:steklov} we may obtain the optimal constant $K_{\dOm,\Om}$ in Proposition~\ref{prop:mvr} as $\frac{|\Om|}{|\dOm|} \sigma^{-1}$, where $\sigma$ denotes the first nontrivial Steklov eigenvalue. In the present example the first Steklov eigenvalue is $coth(1)-tanh(1/2)$, cf.~\cite{escobar} and thus
\begin{equation*}
K_{\dOm,\Om}=\frac{cosh(1)-1}{sinh(1)} (coth(1)-tanh(1/2))^{-1} \approx 0.5431.
\end{equation*}
Inserting this in Proposition \ref{prop:trivialpc} results in $C_\al \le C_{\dOm}\approx 1.3811$. Furthermore we insert the same set of constants in Proposition \ref{prop:mvr}.\\
Using the general results from the previous section instead, we find different constants: We have $h_1(t)=h_2(t)=cosh(t)-coth(1)sinh(t)$ and $H_1(\{\rho_{\dOm}=t\})=2\pi sinh(t)$ from which follows by Lemma \ref{lem:rhol2} (with $t_0\to0$) and Lemma \ref{lem:rho} (with $t_1\to 1$) that
\begin{equation*}
\forall \eps>0\ \exists \varphi:\ |\nabla \varphi|_2^2\le \eps,\ |\Delta \varphi|_2^2\le \frac{1}{2(cosh(1)-1)}+\eps,
\end{equation*}
and
\begin{equation*}
\exists \rho:\ |\nabla \rho|_\infty \le 1,\ |(\Delta \rho)^-|_\infty \le 2.3131.
\end{equation*}
Inserting this in Proposition \ref{prop:k1k2} and Proposition \ref{prop:kdomom} we get
\begin{equation*}
K'_1 =  \frac{sinh^2(1/2)}{sinh^2(1)}\cdot C_\Om \approx 0.0664  ,\ K'_2=0,\ K'_{\dOm,\Om} \approx 0.8981.
\end{equation*}
Inserting this in Proposition \ref{prop:trivialpc} results in $C_\al \le C_{\dOm}$. Furthermore we insert the same set of constants in Proposition \ref{prop:mvr}.\\

\begin{figure}
	\includegraphics[width=\linewidth]{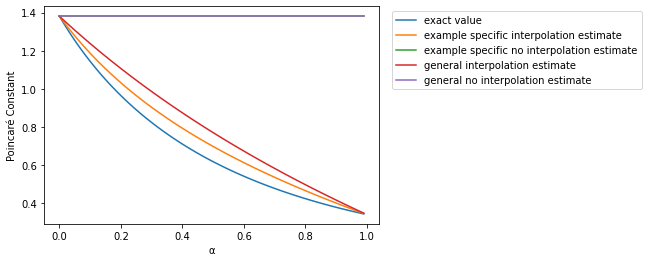}
	\caption{Exact Poincar\'{e} constants (blue), interpolation (yellow) and no interpolation (green) results using computations specifically adapted to the example, interpolation (red) and no interpolation (purple) results using computations not specifically adapted to the example.}
	\label{fig:hpball}
\end{figure}

The curves for actual Poincar\'{e} constants and respective upper bounds via Proposition~\ref{prop:mvr} as well as via Proposition~\ref{prop:trivialpc} obtained by plugging in the quantities collected above are depicted in Figure \ref{fig:hpball}. Note that the green and purple curves overlap.
Again from the figure we may see that our general results are only slightly worse than the ones obtained from computations specifically adapted to the example, and that the interpolation approach results in a significant improvement compared to no interpolation.

\section{Purely sticky reflection case ($\beta=0$)} \label{sec:pinb}
The above results may as well be used to give upper bounds for Brownian motion with sticky reflection from the boundary (but without boundary diffusion). I.e.\ under the same assumptions on $\Om$ as above we consider a diffusion on $\Om$ with Feller generator $(\dD(\hat{A}),\hat{A})$
\begin{align*}
	\dD(\hat{A}) &= \{f\in C(\Om)\ |\ \hat{A}f\in C(\Om)\} \\
	\hat{A}f &= \Del f\1_\Om - \ga \frac{\pd f}{\pd N}\1_{\dOm},
\end{align*}
where $\frac{\pd f}{\pd N}$ is the outer normal derivative and $\ga>0$, which corresponds to inward sticky reflection at $\dOm$. A construction was given again in~\cite{gv} using Dirichlet forms.
We choose $\al\in(0,1)$, such that
\begin{equation*}
	\frac{\alpha}{1-\alpha} \frac{|\dOm|}{|\Om|} = \ga.
\end{equation*}
and set
\begin{equation*}
	\laa := \al \lo + (1-\al)\ld.
\end{equation*}
We find that $-\hat{A}$ is $\la_\al$-symmetric with spectral gap $\hat{\sigma}_\al$ characterised by the Rayleigh quotient resp. Poincar\'{e} constant $\hat{C}_\al$
\begin{equation*}
	\hat{\sigma}_{\al} = \inf_{\substack{ f\in C^1(\Om)\\ Var_{\laa} (f)>0}} \frac{\hat{\eE}_{\al}(f)}{Var_{\laa}(f)}, \qquad \hat{C}_\al := \hat{\sigma}_\al^{-1} = \sup_{\substack{ f\in C^1(\Om)\\ \hat{\eE}_\al(f)>0}} \frac{Var_{\laa}(f)}{\hat{\eE}_\al(f)}
\end{equation*}
where
\begin{equation*}
	Var_{\laa}(f) = \int_\Om f^2 d\laa - \lp\int_\Om f d\laa \rp^2, \hspace{1em} \hat{\eE}_{\al} (f) = \al \int_\Om |\na f|^2 d\lo, \ f\in C^1(\Om)
\end{equation*}
By $C_\Om$ and $C_{\dOm}$ we still denote the usual (Neumann) Poincar\'{e} constants of $\Om$ and $\dOm$ respectively, which we assume to be known.\\

The eigenvalue problem corresponding to the Poincar\'{e} constant may be stated as
\begin{align*}
	\begin{cases}
		\Del f = -\la f  & \text{ in } \Om,\\
		-\gamma \frac{\partial f}{\partial N}=-\la f & \text{ on } \dOm.
	\end{cases}
\end{align*}
This type of eigenvalue problems with eigenvalue featured in the boundary condition has been of separate interest, see e.g.~\cite{binding},~\cite{below},~\cite{gilles},~\cite{barbu},~\cite{buoso}. For Brownian motion with sticky reflection spectral asymptotics have been examined e.g. in~\cite{below}, however we are not aware of results on the spectral gap.\\

\begin{proposition}\label{prop:pcwobd}
Assume there exist constants $K_{\dOm,\Om},K_1$ such that for any $f\in C^1(\Om)$
\begin{equation*}
	Var_{\ld} (f) \le K_{\dOm, \Om} \int_{\Om} |\nabla f|^2 d\lo
\end{equation*}
and
\begin{equation*}
\lp\int_{\Om} f d\lo - \int_{\dOm} f d\ld\rp^2 \le K_1 \int_\Om |\nabla f|^2 d\lo,
\end{equation*}
then it holds for any $\al\in(0,1)$
\begin{equation*}
\hat{C}_\al \le C_\Om + \frac{(1-\al)}{\al} K_{\dOm,\Om} + (1-\al) K_1 .
\end{equation*}
\end{proposition}

\begin{proof}
Let $f\in C^1(\Om)$
\begin{align*}
Var_{\laa}(f) &= \al Var_{\lo}(f) + (1-\al) Var_{\ld}(f) + \al(1-\al) \lp \int_\Om f d\lo - \int_{\dOm} f d\ld \rp^2 \\
	&\le \al C_\Om \int_\Om |\nabla f|^2 d\lo + (1-\al) K_{\dOm,\Om} \int_\Om |\nabla f|^2 d\lo + \al(1-\al) K_1 \int_\Om |\nabla f|^2 d\lo \\
	& = \lp C_\Om + \frac{(1-\al)}{\al} K_{\dOm,\Om} + (1-\al) K_1\rp \al \int_\Om |\nabla f|^2 d\lo.
\end{align*}
\end{proof}

For $f\in C^1(\Om)$ with $Var_{\ld} (f)>0$ the term $Var_{\laa}(f)$ stays positive as $\al$ tends to zero while $\hat{\eE}_{\al}(f)$ vanishes. Thus $\hat{C}_\al$ blows up as $\al$ tends to zero and accordingly so does the bound on $\hat{C}_\al$ proven in Proposition \ref{prop:pcwobd}. Of course the interpolation approach from above is not of any use anymore in this setting. As we have previously shown that inequality~\eqref{eq:mixedterm} can be fulfilled with $K_2=0$, we may now insert $K_1$ and $K_{\dOm,\Om}$ as computed above.
Thus in sum the results obtained above for Brownian motion with sticky reflecting boundary diffusion may also be used for the case without boundary diffusion precisely for the reason that we were able to show that inequality~\eqref{eq:mixedterm} is fulfilled with $K_2=0$. If we make the more strict assumption of respective positive lower bounds $k_R$ on Ricci curvature on $\Om$ and $\gamma$ on second fundamental form on $\dOm$ as in~\cite[section 3.2]{mvr} we may as well insert the set of constants $K_1=\frac{d-1}{dk_R}, K_2=0, K_{\dOm,\Om}=\frac{|\Om|}{|\dOm|} \frac{2}{\gamma}$ obtained there. Note that this is possible as we again have $K_2=0$.\\

\subsection{Examples}
We again consider as examples a unit ball in $\mathbb{R}^2$ and a unit metric ball in the hyperbolic plane.

\subsubsection{Ball in $\mathbb{R}^2$}

\label{ex:euclidean wobd}
As in Example \ref{ex:euclidean} we consider a unit ball in $\mathbb{R}^2$. In order to compute the exact values for $\hat{C}_\al, \al\in(0,1)$ we again proceed as described in \cite[section 3.1]{mvr} and only need to adapt the boundary condition. I.e.\ an eigenfunction $f$ fulfills
\begin{align*}
	\begin{cases}
		\Del f = -\la f  & \text{ in } \Om\\
		- \frac{2\al}{1-\al} \frac{\pd f}{\pd N} =  \Del f & \text{ on } \dOm
	\end{cases}
\end{align*}
and by passing to polar coordinates and seperating variables we obtain a family $\hat{\lambda}_{m,l}, m,l\in\mathbb{N}_{0}$ characterised by
\begin{equation*}
\sqrt{\la} J_m''(\sqrt{\la}) + J_m'(\sqrt{\la}) \frac{1+\al}{1-\al} - J_m(\sqrt{\la})\frac{m}{\la} = 0,
\end{equation*}
where $J_m, m\in\mathbb{N}_0$ are the Bessel functions of the first kind. We then get $\hat{\lambda}_\al = \min_{m,l\in\mathbb{N}_0} \hat{\lambda}_{m,l}$ and $\hat{C}_\al = \frac{1}{\hat{\la}_\al}$.

\begin{figure}
	\includegraphics[width=0.75\linewidth]{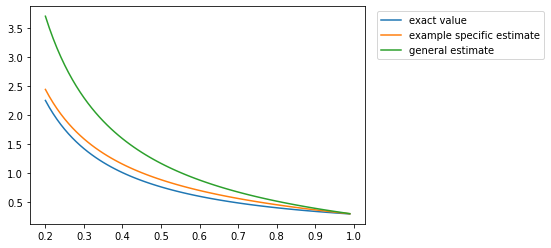}
	\caption{Exact Poincar\'{e} constants (blue), estimates using computations specifically adapted to the example (yellow), estimates using computations not specifically adapted to the example (green).}
	\label{fig:r2ballwobd}
\end{figure}

\smallskip
In order to calculate the explicit values for the upper bound stated in Proposition \ref{prop:pcwobd} we need $C_\Om, K_{\dOm,\Om}$ and $K_1$. All of these have been computed in Example \ref{ex:euclidean}, in particular $K_{\dOm,\Om}$ and $K_1$ have been computed once in a manner adapted to the specific example and once from the previously stated general results (the latter marked by $'$).\\
The curves for actual Poincar\'{e} constants and respective upper bounds via Proposition~\ref{prop:pcwobd} obtained by plugging in the quantities collected above are depicted in Figure \ref{fig:r2ballwobd}. As mentioned before $\hat{C}_\al$ blows up as $\al$ tends to zero. We therefore only consider $\al\ge 0.2$ for the plot.
Figure \ref{fig:r2ballwobd} suggests that Proposition \ref{prop:pcwobd} offers a precise upper bound for $\hat{C}_\al$ that is (in particular for small $\al$) highly depended on how close the values of $K_1$ and $K_{\dOm,\Om}$ are to the optimal constants in inequality~\eqref{eq:varboundary} and inequality~\eqref{eq:mixedterm}. More precisely for small values of $\alpha$ it is mainly the precision of the value for $K_{\dOm,\Om}$ that is relevant. Note that the value $K'_{\dOm,\Om}$ obtained from our general results is worse than $K_{\dOm,\Om}$ obtained from computations adapted to the specific example, while the opposite is true for the values of $K'_1$ and $K_1$.

\subsubsection{Ball in hyperbolic plane}
\label{ex: hyperbolic wobd}
As in Example \ref{ex:hyperbolic} we consider a unit ball in the hyperbolic plane and use the notation introduced above. To calculate the exact values of the Poincar\'{e} constants we proceed as explained in Example \ref{ex:hyperbolic} and only need to adapt the boundary condition (see the second equation in~\eqref{eq:hypsyst}) to
\begin{equation*}
\frac{-\al}{1-\al} \frac{sinh(1)}{cosh(1)-1}\frac{\partial f}{\partial N} = \Delta f \text{ on } \dOm.
\end{equation*}
Inserting $f(t,\xi)= T(t)G(\xi)$ and then using $T(t)= P^\mu_\nu(cosh(t))$ as before, this results in
\begin{gather*}
\frac{-\al}{1-\al} \frac{sinh(1)}{cosh(1)-1} T'(1)G(\xi) = \lp \frac{cosh(1)}{sinh(1)} T'(1)+ T''(1) \rp G(\xi) + \frac{T(1)}{sinh^2(1)}\square G(\xi) \\
\Leftrightarrow T'(1)\lp \frac{\al}{1-\al} \frac{sinh(1)}{cosh(1)-1} + \frac{cosh(1)}{sinh(1)} \rp + T''(1) - \gamma \frac{T(1)}{sinh^2(1)} =0 \Leftrightarrow\\
 {P^\mu_\nu}'' (cosh(1))sinh^2(1) + {P^\mu_\nu}'(cosh(1))\lp \frac{\al}{1-\al}\frac{sinh^2(1)}{cosh(1)-1} + 2cosh(1) \rp - \gamma \frac{P^\mu_\nu(cosh(1))}{sinh^2(1)} = 0.
\end{gather*}
We obtain a two-parameter family of values $\hat{\lambda}_{l,n}$ satisfying this equation. For $\al\in(0,1)$ $\hat{\lambda}_\al := \min_{l,n} \hat{\lambda}_{l,n}$ is the desired spectral gap and $\hat{C}_\al=1/\hat{\lambda}_\al$.

\begin{figure}
	\includegraphics[width=0.75\linewidth]{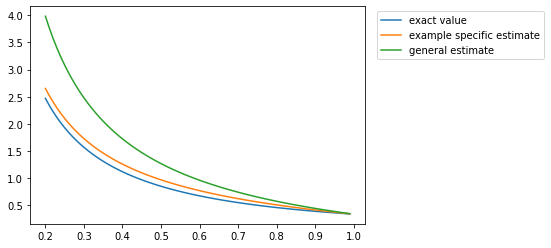}
	\caption{Exact Poincar\'{e} constants (blue), estimates using computations specifically adapted to the example (yellow) and estimates using computations not specifically adapted to the example (green).}
	\label{fig:hpballwobd}
\end{figure}

\smallskip
In order to calculate the explicit values for the upper bound stated in Proposition \ref{prop:pcwobd} we may again use the values for $C_\Om, K_{\dOm,\Om}$ and $K_1$ as computed in Example \ref{ex:hyperbolic} in a manner adapted to the specific example or from the previously stated general results (the latter marked by $'$). The curves for actual Poincar\'{e} constants and respective upper bounds via Proposition~\ref{prop:pcwobd} obtained by plugging in these quantities are depicted in Figure \ref{fig:hpballwobd}. Again we only consider $\al\ge 0.2$ for the plot, as $\hat{C}_\al$ blows up as $\al$ tends to zero.
From Figure \ref{fig:hpballwobd} we may again see that the precision of the upper bound for $\hat{C}_\al$ offered in Proposition \ref{prop:pcwobd} depends particularly for small $\al$ highly on how close the values of $K_1$ and $K_{\dOm,\Om}$ are to the optimal constants in inequality~\eqref{eq:varboundary} and inequality~\eqref{eq:mixedterm}. Note that in this example again the value $K'_{\dOm,\Om}$ is worse than $K_{\dOm,\Om}$, while the opposite is true for $K'_1$ and $K_1$.

\section{Logarithmic Sobolev Inequality}\label{sec:lsi}

Using the notation from above we say that a (tight) logarithmic Sobolev inequality is fulfilled if
\begin{equation*}
\exists L_\al\ge 0 \text{ s.t. } Ent_{\laa}(f^2) \le L_\al \cdot \eE_\al(f)\ \forall f\in C^1(\Om).
\end{equation*}
where
\begin{align*}
	Ent_{\laa}(f^2) &= \int_\Om f^2 log(f^2) d\laa - \lp\int_\Om f^2 d\laa \rp log\lp \int_\Om f^2 d\laa\rp \\
	\eE_{\al} (f) &= \al \int_\Om |\na f|^2 d\lo + (1-\al) \int_{\dOm} |\na^\tau f|^2 d\ld,\ f\in C^1(\Om)
\end{align*}
and $\na^\tau$ denotes the tangential derivative operator on $\dOm$.

\smallskip

In the following by $L_\al$ we will denote the optimal such constant. By $L_\Om$ respectively $L_{\dOm}$ we will denote the optimal logarithmic Sobolev constant associated to the Laplace operator on $\Om$ with Neumann boundary conditions and the logarithmic Sobolev constant associated to the Laplace-Beltrami operator on $\dOm$. We assume that $L_\Om$ and $L_{\dOm}$ (or upper bounds for them) are known. We aim at bounding $L_\al$ for $\al\in(0,1)$ from above.

\smallskip
We consider here the entropy with respect to $\laa$ which is a mixture or more specifically a convex combination of the two measures $\lo$ and $\ld$. The entropy with respect to mixtures of two measures such as $\laa$ has been considered previously e.g. in~\cite{chafai},~\cite{schlichting1},~\cite{schlichting2}.
We first show an analogue of~\cite[Proposition 2.1]{mvr} for the entropy with respect to $\laa$:

\begin{proposition}\label{prop:logsobinterpol}
Assume there are constants $K_{\dOm,\Om}, L_{\dOm,\Om}, K_1, K_2$ such that $\forall f\in C^1(\Om)$:
\begin{align}
Var_{\ld} (f) &\le K_{\dOm,\Om} \int_{\Om} |\nabla f|^2 d\lo \label{eq:var}\\
Ent_{\ld}(f^2) &\le L_{\dOm,\Om} \int_{\Om} |\nabla f|^2 d\lo \label{eq:ent}\\
\lp \int_\Om f d\lo - \int_{\dOm} f d\ld \rp^2 &\le K_1 \int_\Om |\nabla f|^2 d\lo + K_2 \int_{\dOm} |\nabla^\tau f|^2 d\ld \label{eq:mixed}.
\end{align}
Then it holds for any $\al \in(0,1)$
\begin{align*}
L_\al & \le  \inf_{s,t\in[0,1]} \max \Bigl\{  L_\Om + \frac{(1-\al)}{\al}sL_{\dOm,\Om} + \frac{(1-\al)(log(\al)-log(1-\al))}{2\al-1} \lp C_\Om + tK_{\dOm,\Om} + K_1 \rp, \\
	 &\phantom{ \le \inf_{s,t\in[0,1]} \max \Bigl\{ } (1-s)L_{\dOm} + \frac{\al(log(\al)-log(1-\al))}{2\al-1}\lp (1-t)C_{\dOm} + K_2 \rp \Bigr\}.
\end{align*}
\end{proposition}

\begin{proof}
Applying~\eqref{eq:var} and~\eqref{eq:ent} we can estimate for any $f\in C^1(\Om)$:
\begin{align*}
Var_{\ld}(f) \le t K_{\dOm,\Om} \int_{\Om} |\nabla f|^2 d\lo + (1-t)C_{\dOm}  \int_{\dOm} |\nabla^\tau f|^2 d\ld \\
Ent_{\ld}(f^2) \le s L_{\dOm,\Om} \int_{\Om} |\nabla f|^2 d\lo + (1-s)L_{\dOm}  \int_{\dOm} |\nabla^\tau f|^2 d\ld
\end{align*}
for any $s,t\in[0,1].$
We apply this in the following after first using a decomposition of the entropy with respect to the mixture of two measures as well as an optimal logarithmic Sobolev inequality for Bernoulli measures as described in~\cite[section 4]{chafai}
\begin{align*}
E&nt_{\laa}(f^2) \le \al Ent_{\lo}(f^2) + (1-\al) Ent_{\ld}(f^2) \\
	&~+ \frac{\al(1-\al)(log(\al)-log(1-\al))}{2\al-1} \lp Var_{\lo}(f) + Var_{\ld}(f) + (\E_{\lo}(f)-\E_{\ld}(f))^2\rp\\
	 &\le \int_{\Om} |\nabla f|^2 d\lo \lp \al L_\Om + (1-\al)sL_{\dOm,\Om} + \frac{\al(1-\al)(log(\al)-log(1-\al))}{2\al-1} \lp C_\Om + tK_{\dOm,\Om} + K_1 \rp \rp \\
	&+ \int_{\dOm} |\nabla^\tau f|^2 d\ld \lp (1-\al)(1-s)L_{\dOm} + \frac{\al(1-\al)(log(\al)-log(1-\al))}{2\al-1}\lp (1-t)C_{\dOm} + K_2 \rp \rp.
\end{align*}
And thus
\begin{align*}
L_\al & \le  \inf_{s,t\in[0,1]} \max \Bigl\{  L_\Om + \frac{(1-\al)}{\al}sL_{\dOm,\Om} + \frac{(1-\al)(log(\al)-log(1-\al))}{2\al-1} \lp C_\Om + tK_{\dOm,\Om} + K_1 \rp, \\
	 &\phantom{ \le \inf_{s,t\in[0,1]} \max \Bigl\{ } (1-s)L_{\dOm} + \frac{\al(log(\al)-log(1-\al))}{2\al-1}\lp (1-t)C_{\dOm} + K_2 \rp \Bigr\}.
\end{align*}
\end{proof}

As mentioned previously one may again add another optimisation over the pairs of values $K_1, K_2$ fulfilling~\eqref{eq:mixed} in order to improve the upper bound on $L_\al$.

\begin{remark}
Using that for any $a,b,c,d,e,\theta\in\mathbb{R}_{\ge 0}$ it holds
\begin{align*}
&\inf_{s,t\in [0,1]} \max (a+sb+tc, d-se-t\theta)\\
	&=\begin{cases}
		a, &\text{ if } a>d \\
		d-e-\theta, &\text{ if } d-e-\theta>a+b+c \\
		a+c +b\cdot \lp\frac{d-a-c-\theta}{e+b}\rp, &\text{ if } a\le d, d-e-\theta \le a+b+c, b-c\cdot\frac{e+b}{c+\theta}\ge 0, \frac{d-a-c-\theta}{e+b}\ge0 \\
		a+c\cdot\lp\frac{d-a}{c+\theta}\rp, &\text{ if }  a\le d, d-e-\theta \le a+b+c, b-c\cdot\frac{e+b}{c+\theta}\ge 0, \frac{d-a-c-\theta}{e+b} < 0\\
		a + b\cdot\lp \frac{d-a}{e+b} \rp, &\text{ if } a\le d, d-e-\theta \le a+b+c, b-c\cdot\frac{e+b}{c+\theta}< 0, \frac{d-a}{e+b}\le 1 \\
		a+b + c\cdot\lp \frac{d-a-e-b}{c+\theta} \rp, &\text{ if }  a\le d, d-e-\theta \le a+b+c, b-c\cdot\frac{e+b}{c+\theta}< 0, \frac{d-a}{e+b} > 1
	\end{cases}\\
	&= \max\Bigl\{a, d-e-\theta,  \min\bigl[a+c + b \lp\frac{d-a-c-\theta}{e+b}\rp, \\  & \phantom{=\max\Bigl\{a, d-e-\theta,  \min\bigl[}
  a+c \lp\frac{d-a}{c+\theta}\rp, a+ b \lp \frac{d-a}{e+b}\rp, a+b + c \lp \frac{d-a-e-b}{c+\theta} \rp\bigr]\Bigr\},
\end{align*}
we may rewrite the result of Proposition \ref{prop:logsobinterpol} in analogy with Proposition \ref{prop:mvr} as follows: For any $\al \in(0,1)$
\begin{align*}
L_\al \le &\max\Bigl\{a, d-e-\theta, \\ &\hspace{0.5cm}\min\Bigl[a+c + b \lp\frac{d-a-c-\theta}{e+b}\rp, a+c \lp\frac{d-a}{c+\theta}\rp, a+ b \lp \frac{d-a}{e+b}\rp, a+b + c \lp \frac{d-a-e-b}{c+\theta} \rp\Bigr]\Bigr\},
\end{align*}
where
\begin{small}
\begin{alignat*}{2}
&a = L_\Om +\frac{ (1-\al)(log(\al)-log(1-\al))}{2\al-1}(C_\Om + K_1),\hspace{0.5cm} &&b = \frac{1-\al}{\al} L_{\dOm,\Om},\\
&c = \frac{(1-\al)(log(\al)-log(1-\al))}{2\al-1} K_{\dOm,\Om}, &&d = L_{\dOm} + \frac{\al(log(\al)-log(1-\al))}{2\al-1}(C_{\dOm} + K_2) \\
&e = L_{\dOm}, &&\theta =  \frac{\al(log(\al)-log(1-\al))}{2\al-1}C_{\dOm}.
\end{alignat*}
\end{small}
\end{remark}

For comparison we also state the simpler result one can obtain based on~\cite[section 4]{chafai} without interpolation.

\begin{proposition}\label{prop:logsobcm}
Assume there exist constants $K_1, K_2$ such that $\forall f\in C^1(\Om)$:
\begin{equation*}
\lp \int_\Om f d\lo - \int_{\dOm} f d\ld \rp^2 \le K_1 \int_\Om |\nabla f|^2 d\lo + K_2 \int_{\dOm} |\nabla^\tau f|^2 d\ld.
\end{equation*}
Then it holds for any $\al \in(0,1)$
\begin{align*}
L_\al \le  \max\Bigl\{& L_\Om + \frac{(1-\al)(log(\al)-log(1-\al))}{2\al-1}(C_\Om+K_1), \\
& L_{\dOm} + \frac{\al(log(\al)-log(1-\al))}{2\al-1}(C_{\dOm}+K_2)\Bigr\}.
\end{align*}
\end{proposition}

\begin{proof}
As in the proof of Proposition \ref{prop:logsobinterpol} we use a decomposition of the entropy with respect to the mixture of two measures as well as an optimal logarithmic Sobolev inequality for Bernoulli measures as described in~\cite[section 4]{chafai}:
\begin{align*}
Ent_{\laa}(f^2) \le &\al Ent_{\lo}(f^2) + (1-\al) Ent_{\ld}(f^2) + \frac{\al(1-\al)(log(\al)-log(1-\al))}{2\al-1} \\
 & ~~\cdot \lp Var_{\lo}(f) + Var_{\ld}(f) + (\E_{\lo}(f)-\E_{\ld}(f))^2\rp\\
	\le  & \al L_\Om \int_\Om |\na f|^2 d\lo + (1-\al) L_{\dOm}  \int_{\dOm} |\na^\tau f|^2 d\ld + \frac{\al(1-\al)(log(\al)-log(1-\al))}{2\al-1} \\
 & ~~\cdot  \lp (C_\Om + K_1) \int_\Om |\na f|^2 d\lo + (C_{\dOm}+K_2)  \int_{\dOm} |\na^\tau f|^2 d\ld   \rp \\
	= & \lp L_\Om + \frac{(1-\al)(log(\al)-log(1-\al))}{2\al-1} (C_\Om + K_1)  \rp \al \int_\Om |\na f|^2 d\lo \\
	&+ \lp L_{\dOm} + \frac{\al(log(\al)-log(1-\al))}{2\al-1} (C_{\dOm}+K_2)  \rp (1-\al) \int_{\dOm} |\na^\tau f|^2 d\ld   \\
	\le & \max( L_\Om + \frac{(1-\al)(log(\al)-log(1-\al))}{2\al-1} (C_\Om + K_1), \\
	& \phantom{\max (} L_{\dOm} + \frac{\al(log(\al)-log(1-\al))}{2\al-1} (C_{\dOm}+K_2)) \cdot \eE_\al (f).
\end{align*}
\end{proof}

In general  the logarithmic Sobolev constant of a mixture of two measures may blow up as the mixture proportion goes to 0 or 1, c.f.\ \cite{chafai}. Accordingly so may our bounds for the logarithmic Sobolev constant as $\al$ approaches 0 or 1.  The upper bound in Proposition \ref{prop:logsobcm} always blows up as $\al$ approaches 0 or 1. The same is true for the bound in Proposition \ref{prop:logsobinterpol} as $\al$ tends to 0 but not necessarily as $\al$ tends to 1.
On the contrary the Poincar\'{e} constant as well as the upper bound for it in Proposition \ref{prop:mvr} does not blow up as $\al$ tends to 0 or 1.

\subsection{Example: Ball in $\mathbb R^2$} 
\begin{figure}
\includegraphics[width=0.75\linewidth]{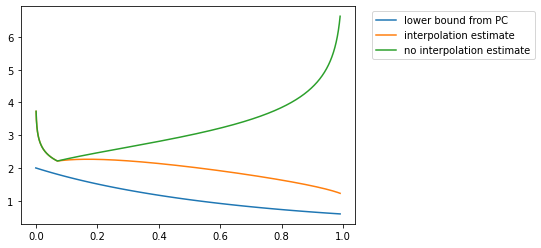}
	\caption{Lower bound via exact Poincar\'{e} constants (blue), and upper bound via interpolation (yellow) and no interpolation (green) for the Logarithmic Sobolev constant }
	\label{fig:r2balllsi}
\end{figure}

To show the feasibility of the approach consider again Example \ref{ex:euclidean}: We consider $\Om=B_1$, the unit ball in $\mathbb{R}^2$. We evaluate the upper bounds for $L_\al, \al \in(0,1)$ obtained via Proposition \ref{prop:logsobcm} and Proposition \ref{prop:logsobinterpol}.
Furthermore we use the following quantities collected in Example \ref{ex:euclidean}:
\begin{equation*}
C_\Om = \frac{1}{3.39}, C_{\dOm}= 1, K_1' = \frac{C_\Om}{4},K_2= 0, K_{\dOm,\Om}= \frac{1}{2}.
\end{equation*}
Added to that we need values for $L_\Om, L_{\dOm}$ and $L_{\dOm,\Om}$.
It follows by the Bakry-\'{E}mery criterion that $L_{\dOm} \le 2$, cf.~\cite{bgl}. Furthermore from~\cite{beckner} we obtain $L_{\dOm,\Om}=1$ and from~\cite{wang3} we obtain $L_\Om \le 1.1799$ \\
Since we have not computed the precise value of $L_\al$ we compare our estimate with the corresponding spectral gap, using that a logarithmic Sobolev inequality with constant $C$ implies a Poincar\'{e} inequality with constant $C/2$, cf.~\cite[Proposition 5.1.3]{bgl}. Thus a lower bound for the optimal logarithmic Sobolev constants is given via
\begin{equation*}
L_\al \ge 2\cdot C_\al, \forall \al \in (0,1),
\end{equation*}
with $C_\al$ as computed in Example \ref{ex:euclidean}.\\
We insert the collected quantities into Proposition \ref{prop:logsobcm} and Proposition \ref{prop:logsobinterpol} and depict the results in Figure \ref{fig:r2balllsi}. Note that the yellow and green curves partly overlap.
The figure shows that the interpolation results clearly differ from the no interpolation results and give significantly better bounds than the approach without interpolation. In particular the interpolation approach gives an upper bound that does not blow up as $\al$ tends to 1.

\subsection{Example: Ball in hyperbolic plane} 
\begin{figure}
\includegraphics[width=0.75\linewidth]{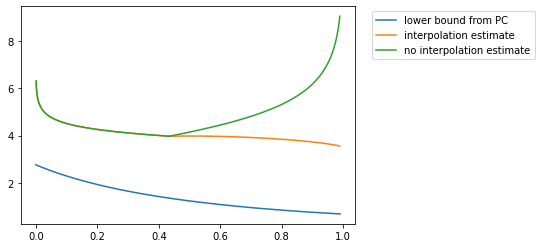}
	\caption{Lower bound via exact Poincar\'{e} constants (blue), and upper bound via interpolation (yellow) and no interpolation (green) for the Logarithmic Sobolev constant }
	\label{fig:hypballlsi}
\end{figure}

We revisit Example \ref{ex:hyperbolic}: We consider the unit metric ball in the hyperbolic plane and evaluate the upper bounds for $L_\al, \al \in(0,1)$ obtained via Proposition \ref{prop:logsobcm} and Proposition \ref{prop:logsobinterpol}. We recall the quantities collected above in Example \ref{ex:hyperbolic}:
 \begin{align*}
C_\Om = 0.3377, C_{\dOm}= sinh^2(1), K_1' =\frac{sinh^2(1/2)}{sinh^2(1)}\cdot C_\Om,\\
 K_2= 0, K_{\dOm,\Om}= \frac{cosh(1)-1}{sinh(1)(coth(1)-tanh(1/2))}.
\end{align*}
Furthermore it follows by the Bakry-\'{E}mery criterion that $L_{\dOm} \le 2\cdot sinh^2(1) $, cf.~\cite{bgl} and from~\cite{wang3} we obtain $L_\Om \le 3.5088$.\\
To obtain an upper bound on $L_{\dOm,\Om}$ we proceed as follows: Rewriting integrals over $\dOm$ as integrals over the Euclidean unit sphere and using the dual-spectral form of the Hardy-Littlewood-Sobolev inequality (for the Euclidean unit sphere) and bound of its right-hand side as presented in~\cite{beckner}[section 2] we obtain for $q\in[2,\infty)$
\begin{equation*}
\lp \int_{\dOm} |f|^q d\ld \rp^{2/q} \le \sum_k a_k^2 (1+(q-2)k) \int_{\dOm} |y_k|^2 d\ld.
\end{equation*}
Here $\sum_k a_k y_k$ is a representation of $f$ in terms of spherical harmonics, i.e. $y_k$ is an eigenfunction of the spherical hyperbolic Laplacian for the eigenvalue $-k^2/sinh^2(1)$.
Using separation of variables we may extend each $y_k$ to a harmonic function on $\Om$ with normal derivative on $\dOm$ equal to $k/sinh(1)$. We thus get as an analogue of~\cite{beckner}[equation (36)] for our unit ball in the hyperbolic plane:
\begin{equation*}
\lp \int_{\dOm} |f|^q d\ld \rp^{2/q} \le \int_{\dOm} |f|^2 d\ld + (q-2)(cosh(1)-1)\int_{\Om} |\nabla u|^2 d\lo,
\end{equation*}
where $u$ is the harmonic extension of $f$ to $\Om$.
Proceeding as in~\cite[Proposition 6.2.3]{bgl} (see also \cite[Proposition 5.1.8]{bgl} for details), this results in
\begin{equation*}
Ent_{\ld}(f^2) \le 2(cosh(1)-1) \int_{\Om} |\nabla f|^2 d\lo,
\end{equation*}
i.e. $L_{\dOm,\Om}\le 2(cosh(1)-1)$.\\
We do not compute the precise value of $L_\al$, but use the lower bound for the optimal logarithmic Sobolev constants given via
\begin{equation*}
L_\al \ge 2\cdot C_\al, \forall \al \in (0,1),
\end{equation*}
with $C_\al$ as computed in Example \ref{ex:hyperbolic}.\\
We insert the collected quantities into Proposition \ref{prop:logsobcm} and Proposition \ref{prop:logsobinterpol} and depict the results in Figure \ref{fig:hypballlsi}. Note that the yellow and green curves partly overlap. Again this figure shows the advantage of the interpolation results as compared to the approach without interpolation.

\section{Purely sticky reflection case ($\beta=0$)}

As in section \ref{sec:pinb} instead of Brownian motion with sticky reflecting boundary diffusion the above results may as well be used to give upper bounds for Brownian motion with sticky reflection from the boundary (but without boundary diffusion).
The Logarithmic Sobolev inequality in this setting is
\begin{equation*}
Ent_{\laa}(f^2) \le \hat{L}_\al \hat{\eE}_\al(f) \ \forall f\in C^1(\Om),
\end{equation*}
where
\begin{equation*}
	\hat{\eE}_{\al} (f) = \al \int_\Om |\na f|^2 d\lo, \ f\in C^1(\Om).
\end{equation*}

\begin{proposition}
Assume there are constants $K_{\dOm,\Om}, L_{\dOm,\Om},K_1$ such that for $f\in C^1(\Om)$
\begin{gather*}
	Var_{\ld} f \le K_{\dOm, \Om} \int_{\Om} |\nabla f|^2 d\lo,\\
	Ent_{\ld}(f^2) \le L_{\dOm,\Om} \int_{\Om} |\nabla f|^2 d\lo,\\
	\lp\int_{\Om} f d\lo - \int_{\dOm} f d\ld\rp^2 \le K_1 \int_\Om |\nabla f|^2 d\lo,
\end{gather*}
then it holds for any $\al\in(0,1)$
\begin{equation*}
\hat{L}_\al \le \lp  L_\Om + \frac{(1-\al)}{\al}L_{\dOm,\Om} + \frac{(1-\al)(log(\al)-log(1-\al))}{2\al-1} \lp C_\Om + K_{\dOm,\Om} + K_1 \rp \rp.
\end{equation*}
\end{proposition}

\begin{proof}
As above we use a decomposition of the entropy with respect to the mixture of two measures as well as an optimal logarithmic Sobolev inequality for Bernoulli measures as described in~\cite[section 4]{chafai}:
\begin{align*}
&\mbox{Ent}_{\laa}(f^2) \le \al Ent_{\lo}(f^2) + (1-\al) Ent_{\ld}(f^2) \\
	&\hspace{1.5cm}+ \frac{\al(1-\al)(log(\al)-log(1-\al))}{2\al-1} \lp Var_{\lo}(f) + Var_{\ld}(f) + (\E_{\lo}(f)-\E_{\ld}(f))^2\rp \\
	&\le \lp  L_\Om + \frac{(1-\al)}{\al}L_{\dOm,\Om} + \frac{(1-\al)(log(\al)-log(1-\al))}{2\al-1} \lp C_\Om + K_{\dOm,\Om} + K_1 \rp \rp \al \int_\Om |\nabla f|^2 d\lo.
\end{align*}
\end{proof}

\section{Boundary-Interior Inequalities}\label{sec:bii}

To obtain explicit bounds on $L_\al$ in the general setting we can use the results from section \ref{sec:pi} for estimate~\eqref{eq:var} and estimate~\eqref{eq:mixed} and so it remains to find $L_{\dOm,\Om}$ such that inequality~\eqref{eq:ent} is fulfilled. In this final section we present a general approach to establish explicit geometric estimates for this quantity under general curvature assumptions.
\subsection{Sobolev-Poincar\'{e}-Trace Inequalities} In the following we present some boundary-interior inequalities that can be proved in a similar fashion as Proposition \ref{prop:k1k2}. They may be seen as alternatives for Proposition \ref{prop:k1k2} but might also be of independent interest. We start from the following statement obtained in the proof of Proposition \ref{prop:kdomom}.

\begin{proposition}\label{prop:boundaryinterior1}
For any $\rho\in C^1(\Om)$ such that $\frac{\partial \rho}{\partial N}\vert_{\dOm}=-1$ and $\nabla \rho$ is Lipschitz continuous on $\Om$ it holds
\begin{align*}
\int_{\dOm} f^2 d\ld &\le \frac{|\Om|}{|\dOm|} \left\{2|\nabla \rho|_\infty \lp \int_\Om f^2 d\lo \cdot \int_\Om |\nabla f|^2 d\lo\rp^{1/2} + |(\Delta \rho)^-|_\infty \int_\Om f^2 d\lo\right\} \\
	&\le \left\{2 |\nabla \rho|_\infty C_\Om^{1/2} + |(\Delta \rho)^-|_\infty C_\Om\right\} \frac{|\Om|}{|\dOm|} \int_\Om |\nabla f|^2 d\lo
\end{align*}
$\forall f\in C^1(\Om)$ with $\int_{\Om} f d\lo =0.$
\end{proposition}

The statement in Proposition \ref{prop:boundaryinterior1} is stronger than needed for inequality~\eqref{eq:mixedterm} in Proposition \ref{prop:mvr} because we bound from above the integral of the squared function as opposed to the square of the integral. Nevertheless the proof of the next corollary follows directly as we may assume for inequality~\eqref{eq:mixedterm} without loss of generality that $\int_\Om f d\lo = 0$. We thus get an upper bound for $K_1$ in inequality~\eqref{eq:mixedterm} that is alternative to Proposition \ref{prop:k1k2}.

\begin{corollary}\label{cor:k1k2}
For any $\rho\in C^1(\Om)$ such that $\frac{\partial \rho}{\partial N}\vert_{\dOm}=-1$ and $\nabla \rho$ is Lipschitz continuous on $\Om$ inequality~\eqref{eq:mixedterm} in Proposition \ref{prop:mvr} is fulfilled with $K_2=0$ and
\begin{equation*}
K_1 = \left\{2 |\nabla \rho|_\infty C_\Om^{1/2} + |(\Delta \rho)^-|_\infty C_\Om\right\}\frac{|\Om|}{|\dOm|}.
\end{equation*}
\end{corollary}

Additionally we think that this computation is of independent interest for the following reason:

\begin{remark}
We may calculate as in the proof of Proposition \ref{prop:boundaryinterior1} to obtain for $f\in C^1(\Om)$ (but not necessarily centered on $\Om$):
\begin{align*}
\int_{\dOm} f^2 d\ld &\le \frac{|\Om|}{|\dOm|} \left\{2|\nabla \rho|_\infty \lp \int_\Om f^2 d\lo \cdot \int_\Om |\nabla f|^2 d\lo\rp^{1/2} + |(\Delta \rho)^-|_\infty \int_\Om f^2 d\lo\right\} \\
	&\le \frac{|\Om|}{|\dOm|} \left\{|\nabla \rho|_\infty  \int_\Om |\nabla f|^2 d\lo  + \lp |\nabla \rho|_\infty  + |(\Delta \rho)^-|_\infty\rp \int_\Om f^2 d\lo\right\}.
\end{align*}
From this it follows that for $K_3:= \frac{|\Om|}{|\dOm|}\lp |\nabla \rho|_\infty  + |(\Delta \rho)^-|_\infty \rp$
\begin{equation*}
|f|_{L^2(\dOm,\ld)}^2 \le K_3 |f|_{W^{1,2}(\Om,\lo)}^2 \Leftrightarrow |f|_{L^2(\dOm,\ld)} \le \sqrt{K_3} |f|_{W^{1,2}(\Om,\lo)}.
\end{equation*}
As $W^{1,2}(\Om,\lo)$ is the completion of smooth functions whose derivatives up to degree 1 are in $L^2(\Om,\lo)$, the inequality also holds for all functions in $W^{1,2}(\Om,\lo).$
Thus via stating a specific constant $K_3$, as can be obtained from Lemma \ref{lem:rho}, we also give an upper bound for the norm of the Trace operator $\vert_{\dOm}:W^{1,2}(\Om,\lo)\to L^2(\dOm,\ld)$ that is explicit in terms of upper bounds on sectional curvature and second fundamental form on the boundary. (An optimal upper bound in terms of the geometry of $\Om$ seems to be unknown in this form as of yet.)
\end{remark}

\begin{proposition}\label{prop:bdrytraceopest}
Let $k_2\in\mathbb{R}$ such that $sect \le k_2$ and $\gamma_2\in\mathbb{R}$ such that $\Pi\le \gamma_2 id$. Then the norm of the Trace operator $\vert_{\dOm}:W^{1,2}(\Om,\lo)\to L^2(\dOm,\ld)$ is bounded from above by
\begin{align*}
\lp \frac{|\Om|}{|\dOm|}\lp 1 + \inf_{t_1\in(0,h_2^{-1}(0))} \sup_{t\in(0,t_1)} \lp(d-1)\frac{h_2'}{h_2}\lp t\rp \lp 1-\frac{t}{t_1}\rp - \frac{1}{t_1}\rp^-\rp\rp^{1/2}.
\end{align*}
\end{proposition}

It is known that on a smooth, compact $d$-dimensional Riemannian manifold $(\Om,g)$ for $q\in[1,d)$ and $\frac{1}{p}=\frac{1}{q} - \frac{1}{d}$ (and thus for all $p\in[1,\frac{qd}{d-q}]$)  there is a constant $C_{p,q}$ such that for $f\in H^{1,q}(\Om)$:
\[\lp \int_\Om |f-\bar{f}|^p d\lo \rp^{1/p} \le C_{p,q} \lp \int_\Om |\nabla f|^q d\lo \rp^{1/q},
\]
where $\bar{f}=\int_\Om f d\lo.$

In terms of these Sobolev-Poincar\'{e} constants we may also show a generalisation of Proposition \ref{prop:boundaryinterior1}:
\begin{proposition}\label{prop:boundaryinterior2}
Let $(\Om,g)$ be a smooth, compact Riemannian manifold of dimension $d\ge 3$, with a connected boundary.
For any $\rho\in C^1(\Om)$ such that $\frac{\partial \rho}{\partial N}\vert_{\dOm}=-1$ and $\nabla \rho$ is Lipschitz continuous on $\Om$ it holds $\forall f\in C^1(\Om)$ with $\int_{\Om} f d\lo =0$ and $p\in [2,\frac{2d-2}{d-2}]$
\begin{equation*}
\lp \int_{\dOm} |f|^p d\ld \rp^{2/p} \le \lp \lp \frac{|\Om|}{|\dOm|} |\nabla \rho|_\infty p\rp^{2/p}  C_{2(p-1),2}^{2(p-1)/p}+ \lp\frac{|\Om|}{|\dOm|} |(\Delta \rho)^-|_\infty\rp^{2/p}   C_{p,2}^2 \rp \int_\Om |\nabla f|^2 d\lo.
\end{equation*}
\end{proposition}

\begin{proof}
We may calculate as in the previous proofs to obtain
\begin{align*}
 \Bigl(& \int_{\dOm} |f|^p d\ld \Bigr)^{2/p}  \le \lp \frac{|\Om|}{|\dOm|} |\nabla \rho|_\infty p \int_\Om |f|^{p-1} |\nabla f| d\lo +\frac{|\Om|}{|\dOm|} |(\Delta \rho)^-|_\infty \int_\Om |f|^p d\lo  \rp^{2/p} \\
	\le &\lp \frac{|\Om|}{|\dOm|} |\nabla \rho|_\infty p \lp \int_\Om |f|^{2(p-1)} d\lo \rp ^{1/2} \lp \int_\Om |\nabla f|^2 d\lo \rp ^{1/2} +\frac{|\Om|}{|\dOm|} |(\Delta \rho)^-|_\infty \int_\Om |f|^p d\lo  \rp^{2/p} \\
	\le &\biggl( \frac{|\Om|}{|\dOm|} |\nabla \rho|_\infty p \lp C_{2(p-1),2} \lp \int_\Om |\nabla f|^2 d\lo\rp^{1/2} \rp^{(p-1)} \lp \int_\Om |\nabla f|^2 d\lo \rp^{1/2} \\
	&+ \frac{|\Om|}{|\dOm|} |(\Delta \rho)^-|_\infty \lp C_{p,2} \lp \int_\Om |\nabla f|^2 d\lo\rp^{1/2} \rp^{p} \biggr)^{2/p} \\
	\le  &\lp \frac{|\Om|}{|\dOm|} |\nabla \rho|_\infty p \lp C_{2(p-1),2}\lp  \int_\Om |\nabla f|^2 d\lo\rp^{1/2} \rp^{(p-1)} \lp \int_\Om |\nabla f|^2 d\lo \rp ^{1/2}\rp^{2/p} \\
	& + \lp\frac{|\Om|}{|\dOm|} |(\Delta \rho)^-|_\infty\rp^{2/p}   C_{p,2}^2 \int_\Om |\nabla f|^2 d\lo \\
	= &\lp \frac{|\Om|}{|\dOm|} |\nabla \rho|_\infty p\rp^{2/p} \lp C_{2(p-1),2}^2 \int_\Om |\nabla f|^2 d\lo \rp^{(p-1)/p} \lp \int_\Om |\nabla f|^2 d\lo \rp ^{1/p}  \\
	&+ \lp\frac{|\Om|}{|\dOm|} |(\Delta \rho)^-|_\infty\rp^{2/p}   C_{p,2}^2 \int_\Om |\nabla f|^2 d\lo \\
	= &\lp \frac{|\Om|}{|\dOm|} |\nabla \rho|_\infty p\rp^{2/p}  C_{2(p-1),2}^{2(p-1)/p} \int_\Om |\nabla f|^2 d\lo  + \lp\frac{|\Om|}{|\dOm|} |(\Delta \rho)^-|_\infty\rp^{2/p}   C_{p,2}^2 \int_\Om |\nabla f|^2 d\lo \\
	= &\lp \lp \frac{|\Om|}{|\dOm|} |\nabla \rho|_\infty p\rp^{2/p}  C_{2(p-1),2}^{2(p-1)/p}+ \lp\frac{|\Om|}{|\dOm|} |(\Delta \rho)^-|_\infty\rp^{2/p}   C_{p,2}^2 \rp \int_\Om |\nabla f|^2 d\lo.
\end{align*}
Here we have used the Sobolev-Poincar\'{e} inequalities associated with $C_{2(p-1),2}$ and $C_{p,2}$. Note therefor that for $p\in \Bigl[2,\frac{2d-2}{d-2}\Bigr]$ it holds $p, 2(p-1) \in \Bigl[1,\frac{2d}{d-2}\Bigr]$.
\end{proof}

Again explicit constants may be obtained from Lemma \ref{lem:rho} in terms of upper bounds on sectional curvature and second fundamental form on the boundary.

\subsection{Boundary Trace Logarithmic Sobolev Inequalities}

 As point of departure we recall the following lemma, cf.~\cite{rothaus}.

\begin{lemma}[Rothaus' Lemma]\label{lem:rothaus}
Let $f:\dOm\to\mathbb{R}$ be measurable and assume that \\ $\int_{\dOm} f^2 log(1+f^2)d\ld <\infty$. For every $a\in\mathbb{R}$
\begin{equation*}
Ent_{\ld}\lp(f+a)^2\rp \le Ent_{\ld} \lp f^2 \rp + 2\int_{\dOm} f^2 d\ld.
\end{equation*}
\end{lemma}

\begin{lemma}\label{lem:bspent}
If $f\in C^1(\Om)$ fulfills $\int_\Om f d\lo=0$ and if there are constants $\tilde{C}_{p,2}$ such that
\begin{equation}\label{eq:cdn1}
\lp \int_{\dOm} |f|^p d\ld \rp^{2/p} \le \tilde{C}_{p,2} \int_{\Om} |\nabla f|^2 d\lo, \forall p\in\left[2, \frac{2d-2}{d-2}\right],
\end{equation}
then it holds
\begin{equation*}
Ent_{\ld}\lp f^2\rp \le \inf_{p\in\left[2,\frac{2d-2}{d-2} \right]} \frac{p}{p-2} \frac{\tilde{C}_{p,2}}{e} \int_{\Om} |\nabla f|^2 d\lo.
\end{equation*}
\end{lemma}

The proof of this Lemma is adapted from ~\cite[Proposition 6.2.3]{bgl}, see also \cite[Proposition 5.1.8]{bgl} for details.

\begin{proof}
Without loss of generality we may assume $\int_{\dOm} f^2 d\ld=1$ and define
\begin{equation*}
\phi:(0,1]\to\mathbb{R},\ \phi(r):= log\lp \lp \int_{\dOm} |f|^{1/r}d\ld\rp^{r} \rp.
\end{equation*}
$\phi$ is convex and $\phi'\lp \frac{1}{2} \rp= -Ent_{\ld}(f^2)$. Now for $p\in \left[ 2, \frac{2d-2}{d-2}\right]$ via the convexity of $\phi$
\begin{align*}
d\lp \phi\lp\frac{1}{2}\rp - \phi\lp\frac{1}{p}\rp\rp &= d \int_{1/p}^{1/2} \phi'(s)ds \le d \phi'\lp\frac{1}{2}\rp \lp \frac{1}{2}- \frac{1}{p}\rp  \\
\Leftrightarrow   & -Ent_{\ld}(f^2) \ge \frac{2p}{p-2} \lp \phi\lp\frac{1}{2}\rp-\phi\lp \frac{1}{p}\rp \rp \\
\Leftrightarrow &Ent_{\ld}(f^2) \le \frac{p}{p-2}log\lp \lp \int_{\dOm} |f|^p d\ld\rp^{2/p} \rp.
\end{align*}
Inserting inequality~\eqref{eq:cdn1} we obtain
\begin{equation*}
Ent_{\ld}(f^2) \le \frac{p}{p-2} log\lp  \tilde{C}_{p,2} \int_{\Om} |\nabla f|^2 d\lo \rp.
\end{equation*}
We define $\tilde{\phi}:(0,\infty)\to\mathbb{R}, \tilde{\phi}(r):=\frac{p}{p-2}log(\tilde{C}_{p,2}r)$. $\tilde{\phi}$ is concave and we may thus compute
\begin{align*}
Ent_{\ld}(f^2) &\le \tilde{\phi}\lp \int_\Om |\nabla f|^2 d\lo \rp \le \tilde{\phi}(r) + \tilde{\phi}'(r)\lp \int_\Om |\nabla f|^2 d\lo -r\rp \\
& = \tilde{\phi}'(r) \int_\Om |\nabla f|^2 d\lo + \lp\tilde{\phi}(r)-r\tilde{\phi}'(r)\rp.
\end{align*}
Choosing $r=\frac{e}{\tilde{C}_{p,2}}$ the last term vanishes and we obtain
\begin{equation*}
Ent_{\ld}(f^2) \le \tilde{\phi}'\lp \frac{e}{\tilde{C}_{p,2}} \rp \int_\Om |\nabla f|^2 d\lo = \frac{p}{p-2} \frac{\tilde{C}_{p,2}}{e} \int_\Om|\nabla f|^2 d\lo.
\end{equation*}
\end{proof}

\begin{proposition}
Assume that $d\ge3$. For any $\rho\in C^1(\Om)$ such that $\frac{\partial \rho}{\partial N}\vert_{\dOm}=-1$ and $\nabla \rho$ is Lipschitz continuous on $\Om$ inequality~\eqref{eq:ent} in Proposition \ref{prop:logsobinterpol} is fulfilled with
\begin{align*}
L_{\dOm,\Om} = & \inf_{p\in\left[2,\frac{2d-2}{d-2} \right]} \frac{p}{p-2} \frac{1}{e} \lp \lp \frac{|\Om|}{|\dOm|} |\nabla \rho|_\infty p\rp^{2/p}  C_{2(p-1),2}^{2(p-1)/p}+ \lp\frac{|\Om|}{|\dOm|} |(\Delta \rho)^-|_\infty\rp^{2/p}   C_{p,2}^2 \rp \\
	&+ 2\left\{2 |\nabla \rho|_\infty C_\Om^{1/2} + |(\Delta \rho)^-|_\infty C_\Om\right\} \frac{|\Om|}{|\dOm|},\notag
\end{align*}
where $(\cdot)^-$ denotes the negative part of a function.
\end{proposition}

\begin{proof}
Let $f\in C^1(\Om)$ then for $a:=\int_\Om f d\lo$ we define $\tilde{f}:=f-a$ and by Lemma \ref{lem:rothaus} it holds
\begin{equation*}
Ent_{\ld}(f^2) = Ent_{\ld}((\tilde{f}+a)^2) \le Ent_{\ld}(\tilde{f}^2) + 2\int_{\dOm} \tilde{f}^2 d\ld.
\end{equation*}
As $\tilde{f}$ is centered on $\Om$ the assumptions of Lemma \ref{lem:bspent} are fulfilled due to Proposition \ref{prop:boundaryinterior2} and we obtain
\begin{align*}
Ent_{\ld}(\tilde{f}^2)  \le &\inf_{p\in\left[2,\frac{2d-2}{d-2} \right]} \frac{p}{p-2} \frac{1}{e}
\lp \lp \frac{|\Om|}{|\dOm|} |\nabla \rho|_\infty p\rp^{2/p}  C_{2(p-1),2}^{2(p-1)/p}+ \lp\frac{|\Om|}{|\dOm|} |(\Delta \rho)^-|_\infty\rp^{2/p}   C_{p,2}^2 \rp \\
	&\cdot \int_{\Om} |\nabla \tilde{f}|^2 d\lo.
\end{align*}
Furthermore by Proposition \ref{prop:boundaryinterior1}
\begin{equation*}
\int_{\dOm} \tilde{f}^2 d\ld \le \left\{2 |\nabla \rho|_\infty C_\Om^{1/2} + |(\Delta \rho)^-|_\infty C_\Om\right\} \frac{|\Om|}{|\dOm|} \int_\Om |\nabla \tilde{f}|^2 d\lo.
\end{equation*}
Thus we have
\begin{align*}
Ent_{\ld}(f^2)  \le \Biggl( &\inf_{p\in\left[2,\frac{2d-2}{d-2} \right]} \frac{p}{p-2} \frac{1}{e}  \lp \lp \frac{|\Om|}{|\dOm|} |\nabla \rho|_\infty p\rp^{2/p}  C_{2(p-1),2}^{2(p-1)/p}+ \lp\frac{|\Om|}{|\dOm|} |(\Delta \rho)^-|_\infty\rp^{2/p}   C_{p,2}^2 \rp \\
	+ &2\left\{2 |\nabla \rho|_\infty C_\Om^{1/2} + |(\Delta \rho)^-|_\infty C_\Om\right\} \frac{|\Om|}{|\dOm|}  \Biggr) \int_{\Om} |\nabla \tilde{f}|^2 d\lo \\
	=  \Biggl( &\inf_{p\in\left[2,\frac{2d-2}{d-2} \right]} \frac{p}{p-2} \frac{1}{e} \lp \lp \frac{|\Om|}{|\dOm|} |\nabla \rho|_\infty p\rp^{2/p}  C_{2(p-1),2}^{2(p-1)/p}+ \lp\frac{|\Om|}{|\dOm|} |(\Delta \rho)^-|_\infty\rp^{2/p}   C_{p,2}^2 \rp \\
	+ &2\left\{2 |\nabla \rho|_\infty C_\Om^{1/2} + |(\Delta \rho)^-|_\infty C_\Om\right\} \frac{|\Om|}{|\dOm|} \Biggr) \int_{\Om} |\nabla f|^2 d\lo.
\end{align*}
\end{proof}

Combining this with Lemma \ref{lem:rho} again results in an explicit upper bound for $L_{\dOm,\Om}$.

\begin{proposition}\label{prop:bdrlogsob}
Assume that $d\ge3$. Let $k_2\in\mathbb{R}$ such that $sect \le k_2$ and $\gamma_2\in\mathbb{R}$ such that $\Pi\le \gamma_2 id$. Then inequality~\eqref{eq:ent} in Proposition \ref{prop:logsobinterpol} is fulfilled with
\begin{align*}
L_{\dOm,\Om} =& \inf_{p\in\left[2,\frac{2d-2}{d-2} \right]} \frac{p}{p-2} \frac{1}{e} \biggl[ \lp \frac{|\Om|}{|\dOm|} p\rp^{2/p}  C_{2(p-1),2}^{2(p-1)/p} \\
	&+ \lp\frac{|\Om|}{|\dOm|}  \inf_{t_1\in(0,h_2^{-1}(0))} \sup_{t\in(0,t_1)} \lp(d-1)\frac{h_2'}{h_2}\lp t\rp \lp 1-\frac{t}{t_1}\rp - \frac{1}{t_1}\rp^-\rp^{2/p}   C_{p,2}^2 \biggr] \\
	&+ 2\left\{2 C_\Om^{1/2} +  \inf_{t_1\in(0,h_2^{-1}(0))} \sup_{t\in(0,t_1)} \lp(d-1)\frac{h_2'}{h_2}\lp t\rp \lp 1-\frac{t}{t_1}\rp - \frac{1}{t_1}\rp^- C_\Om\right\} \frac{|\Om|}{|\dOm|}.
\end{align*}
\end{proposition}

\subsection{Semigroup Approach}

Lastly we present an alternative approach via semigroup theory to bound $L_{\dOm,\Om}$ and $L_\al$:\\
Let $P_t^N$ be the Neumann semigroup on $\Om$. Then
\begin{equation} \label{0}
|P_t^N|_{1\to\infty}\le 1+ c t^{-\ff d 2},\ \ t>0
\end{equation}
holds for some constant $c>0$, where $|\cdot|_{1\to\infty}$ is the operator norm from $L^1(\la_\Om)$ to $L^\infty(\la_\Om).$ Next, let $\scr D$ be the set of positive functions $1\le \phi\in C_b^2(\Om)$  such that
$$\Pi\ge \ff{\partial log(\phi)}{\partial N}\Big|_{\partial M}.$$
See Lemma 3.5.8 in \cite{W14} for concrete examples of $\phi$. Let
$$K_\phi:= \inf\big\{\Ric(X,X)-(\phi \Del \phi^{-1})(x): x\in \Om, X\in T_x M, |X|=1\big\}.$$
Let $L_{\dOm,\Om}'$ be the smallest constant such that $\forall f\in C^1(\Om)$
$$\int_{\dOm} \lp f^2 log(f^2) +1-f^2 \rp d\ld \le L_{\dOm,\Om}' \int_\Om |\na f|^2 d\lo.$$

\begin{proposition}

\begin{enumerate}
\item[$(1)$] We have $L_{\dOm,\Om}\le L_{\dOm,\Om}'$ and
$$ L_\al \le \frac{\al L_\Om+ (1-\al)L_{\dOm,\Om}'}{\al}.$$
 \item[$(2)$] Let $\phi\in \scr D$ and $\rr\in C_b^2(\Om)$ with $\ff{\partial \rr}{\partial N}\big|_{\partial M}=1.$ If $K_\phi>0$, then
$$\ff{|\dOm|}{|\Om|}L'_{\dOm,\Om}\le L_\Om|(\Del\rr)^+|_\infty+ \sqrt{32}|\na\rr|_\infty |\phi|_\infty^3 \bigg(\int_0^\infty e^{-2K_\phi t}(1+ |log(|P_t^N|_{1\to\infty})|) dt\bigg)^{1/2}.$$
In particular, when $Ric\ge K$ for some constant $K>0$ and $\Pi\ge 0$, the estimate holds for $|\phi|_\infty=1$ and $K_\phi=K.$
\item[$(3)$] In general,
\begin{align*} \ff{|\dOm|}{|\Om|}&L'_{\dOm,\Om}\le L_\Om|(\Del\rr)^+|_\infty\\
&+ \sqrt{32}|\na\rr|_\infty |\phi|_\infty^3 \bigg(\int_0^1 e^{-2K_\phi t}(1+ |log(|P_t^N|_{1\to\infty})|) dt+\int_1^\infty e^{-2K_\phi} e^{-(t-1)/C_\Om} |P_1^N|_{1\to\infty} dt \bigg)^{1/2}.\end{align*}
\end{enumerate}
\end{proposition}

\begin{proof} (1)
By Remark 1.20 in \cite{DS}, we have
\begin{equation}
\label{A}Ent_{\la_{\dOm}}(f^2)= \inf_{t>0} \int_{\dOm } \lp f^2 log (f^2)-f^2 log(t)+t-f^2 \rp d\ld.
\end{equation}
By taking $t=1$ we obtain
$$Ent_{\la_{\dOm}}(f^2)\le \int_{\dOm} \lp f^2 log(f^2)+1-f^2\rp d\ld\le L_{\dOm,\Om}' \int_\Om |\nabla f|^2 d\lo,$$
so that $L_{\dOm,\Om}\le L_{\dOm,\Om}'$.\\
Applying \eqref{A} to $\la_\al$ replacing $\la_{\dOm},$ and assuming that $t=\int_\Om f^2 d\lo=1$,  we obtain
\begin{align*} Ent_{\la_\al}(f^2) &\le \al Ent_{\la_\Om}(f^2)+ (1-\al)\int_{\dOm} \lp f^2 log(f^2)+1-f^2 \rp d\ld\\
&\le \big(\al L_\Om+ (1-\al) L'_{\dOm,\Om}\big)\int_\Om |\na f|^2d\lo.\end{align*}

(2) Letting $\int_\Om f^2 d\lo=1$  and noting that
$$f^2 log(f^2)+1-f^2\ge 0,$$ and using integration by parts, we obtain
\begin{align}\label{1}
\frac{|\dOm|}{|\Om|}\int_{\dOm} &\lp f^2 log(f^2)+1-f^2\rp d\ld = \int_\Om (f^2 log(f^2)+1-f^2)\Del\rr+ log(f^2)\na f^2 \cdot \na\rr d\lo\notag \\
	&\le |(\Del\rr)^+|_\infty Ent_{\la_\Om}(f^2) + 2|\nabla\rr|_\infty  \lp \int_\Om |\na f|^2 d\lo \int_\Om f^2(log(f^2))^2d\lo \rp^{1/2}.
\end{align}
Moreover, since $|P_t^Nf^2 -\int_\Om f^2 d\lo|_2 \to 0$ as $t\to\infty$,
we have
\begin{align}\label{2}
\int_\Om f^2(log (f^2))^2 d\lo&= -\int_0^\infty\Big[\frac{d}{dt} \int_\Om (P_t^Nf^2)(log(P_t^Nf^2))^2 d\lo \Big]  dt \notag\\
	&= 2\int_0^\infty\int_\Om \ff{|\na P_t^N f^2|^2}{P_t^Nf^2} \big(1+log(P_t^N f^2)\big) d\lo dt.
\end{align}
By Theorem 3.6.1(2) in \cite{W14}, for the reflecting diffusion process $X_t$ on $\Om$ generated by $\Del$, and a martingale $M_t$ with quadratic variation $\langle M\rangle_t=\int_0^t |\nabla \log \phi(X_s)|^2 d s,$ we have 
\begin{align*}&|\na P_t^N f^2|^2\le 4|\phi|_\infty^2 \big(\E [|f\nabla f|(X_t) e^{-K_\phi t+\sqrt 2 M_t-2\langle M\rangle_t}]\big)^2\\
&\le4|\phi|_\infty^2 e^{-2K_\phi t}  (P_t^Nf^2) \E[ |\nabla f|^2(X_t)e^{2\sqrt 2 M_t-4\langle M\rangle_t}]
= 4|\phi|_\infty^2 e^{-2K_\phi t}  (P_t^Nf^2) \bar P_t^N |\nabla f|^2,\end{align*}
where, by Girsanov's theorem, $\bar P_t^N$ is the Neumann semigroup on $\Om$ generated by $$\Del+4\nabla\log \phi,$$ which  is symmetric in $L^2(\phi^4 \lo)$ so that
$$\int_\Om \bar P_t^N |\nabla f|^2 d\lo \le \int_\Om \bar P_t^N |\nabla f|^2 \phi^4 d\lo
=\int_\Om   |\nabla f|^2 \phi^4d\lo\le |\phi|_\infty^4 \int_\Om  |\nabla f|^2 d\lo.$$
Thus, 
\begin{equation} \label{3}
\int_\Om \ff{|\na P_t^N f^2|^2}{P_t^Nf^2} \big(1+log (P_t^Nf^2)\big) d\lo \le 4 |\phi|_\infty^6 e^{-2K_\phi t}(1+|log (|P_t^N|_{1\to\infty})|)\int_\Om |\nabla f|^2 d\lo.\end{equation}
Combining this with \eqref{1} and \eqref{2} we obtain the desired estimate. When $\Pi\ge 0$ and $Ric\ge K>0$, we may take $\phi\equiv 1$ so that the estimate holds for $|\phi|_\infty=1$ and $K_\phi=K.$

(3) In general,  by the semigroup property and Poincar\'e inequality, when $t\ge 1$ we have
\begin{align*}
&\int_\Om |\na P_t^N f^2|^2 d\lo \le e^{-(t-1)/C_\Om} \int_\Om |\na P_1^N f^2|^2 d\lo\\
&\le 4|\phi|_\infty^2 e^{-2K_\phi} e^{-(t-1)/C_\Om} \int_\Om (\bar{P}_1^N|\na f|^2) (P_1^N f^2) d\lo\\
&\le 4|\phi|_\infty^6  e^{-2K_\phi} e^{-(t-1)/C_\Om} |P_1^N|_{1\to\infty} \int_\Om |\nabla f|^2d\lo. \end{align*}
Combining this with \eqref{3} and the fact that $\ff{1+log r}r\le 1$ for $r>0$, we obtain
$$\int_\Om \ff{|\na P_t^N f^2|^2}{P_t^Nf^2} \big(1+log(P_t^Nf^2)\big)d \lo \le 4|\phi|_\infty^6 \int_\Om|\na f|^2 d\lo \begin{cases}  e^{-2K_\phi} e^{-(t-1)/C_\Om} |P_1^N|_{1\to\infty},\ \ t\ge 1\\
e^{-2K_\phi t} \big(1+log(|P_t^N|_{1\to\infty}) \big),\ \ t\in (0,1).\end{cases}$$
Therefore, the desired estimate follows from \eqref{1} and \eqref{2}.
\end{proof}

\nocite{MR2448584,MR3951758,zbMATH07168051,zbMATH06591437,kolesnikov,colbgir,zbMATH02183024,zbMATH05215538,MR2029716,Nazarov,Matculevich_Repin,MR2255233,MR1443055,MR121855,MR1971310,reilly,MR4065110}

\bibliographystyle{abbrv}
\bibliography{bib}
\end{document}